\documentclass[11pt,reqno]{amsart}
\usepackage{fullpage,times,graphicx,amssymb,amsmath,psfrag,xcolor,pdfsync,array}
\usepackage{bm}
\newtheorem{theorem}{Theorem}[section]
\newtheorem{proposition}[theorem]{Proposition}

\newtheorem{lemma}[theorem]{Lemma}
\newtheorem{corollary}[theorem]{Corollary}
\renewenvironment{proof}{\textbf{Proof.}}{\QED\bigskip}

\usepackage{algorithmic,algorithm}

\usepackage[square]{natbib}

\definecolor{ddarkbrown}{rgb}{0.5,0.2,0.05} \definecolor{bbluegray}{rgb}{0.05,0,0.5}

\newcommand{\BEAS}{\begin{eqnarray*}}
\newcommand{\EEAS}{\end{eqnarray*}}
\newcommand{\BEA}{\begin{eqnarray}}
\newcommand{\EEA}{\end{eqnarray}}
\newcommand{\BEQ}{\begin{equation}}
\newcommand{\EEQ}{\end{equation}}
\newcommand{\BIT}{\begin{itemize}}
\newcommand{\EIT}{\end{itemize}}
\newcommand{\BNUM}{\begin{enumerate}}
\newcommand{\ENUM}{\end{enumerate}}

\newcommand{\BA}{\begin{array}}
\newcommand{\EA}{\end{array}}

\newcommand{\ones}{\mathbf 1}

\newcommand{\reals}{{\mathbb R}}

\newcommand{\symm}{{\mbox{\bf S}}}  % symmetric matrices
\newcommand{\Tr}{\mathop{\bf Tr}}
\newcommand{\diag}{\mathop{\bf diag}}
\newcommand{\lambdamax}{{\lambda_{\rm max}}}

\newcommand{\idm}{\mathbf{I}}

\newcommand{\Expect}{\mathop{\bf E}}

\newcommand{\QED}{~~\rule[-1pt]{6pt}{6pt}}

\newcommand{\argmin}{\mathop{\rm argmin}}

\newcommand{\var}{\mathop{\bf var}}

\newcommand{\argmax}{\mathop{\rm argmax}}

\usepackage[colorlinks,citecolor=bbluegray,linkcolor=ddarkbrown,urlcolor=blue,breaklinks]{hyperref}
\usepackage[off]{auto-pst-pdf} % If no figure recompilation necessary
\newcommand{\norm}[1]{\|{#1}\|}

\newcommand{\tendsto}{\rightarrow}
\newcommand{\symMatrices}{\symm_n}
\newcommand{\trsp}{^T}
\newcommand{\trace}{\Tr}
\newcommand{\eps}{\epsilon}
\newcommand{\Exp}[1]{\Expect\left[{#1}\right]}
\newcommand{\gO}{O}
\newcommand{\lo}{o}
\newcommand{\id}{\idm}
\newcommand{\mysigma}{\eps}
\newcommand{\mysigmasquared}{\eps}

\newcommand{\complexNormal}[2]{{\mathcal N}_{\mathbb{C}}\left(#1,#2\right)}

\newcommand{\centeredNormal}[1]{{\mathcal N}\left(0,#1\right)}
\newcommand{\goodClass}{{\mathcal G}}
\newcommand{\weakCV}{\Longrightarrow}
\newcommand{\twc}{\mathrm{TW}_2}
\newcommand{\tw}{\mathrm{TW}}
\newcommand{\alphaTilde}{\widetilde{\alpha}}
\newcommand{\betaTilde}{\widetilde{\beta}}
\newcommand{\smoothedF}{\mathsf{F}}

\newcommand{\smallQuantity}{\eta}
\newcommand{\canonicalRandomVector}{z}
\newcommand{\canonicalDetVector}{v}
\newcommand{\coordCanonRandomVector}{\mathsf{z}}
\newcommand{\coordCanonDetVector}{\mathsf{v}}
\newcommand{\orthoMat}{\mathcal{O}}
\newcommand{\localLipConstant}[1]{L\left[#1\right]}
\newcommand{\card}[1]{\textrm{card}\left\{#1\right\}}
\newcommand{\iid}{\stackrel{i.i.d}{\sim}}

\newcommand{\equalInLaw}{\stackrel{\mathcal{L}}{=}}
\begin{document}
\title{A Stochastic Smoothing Algorithm for Semidefinite Programming}

\author{Alexandre d'Aspremont}
\address{CNRS \& D.I., UMR 8548, \vskip 0ex
\'Ecole Normale Sup\'erieure, Paris, France.}
\email{aspremon@ens.fr}

\author{Noureddine El Karoui}
\address{Statistics, U.C. Berkeley. Berkeley, CA 94720.}
\email{nkaroui@stat.berkeley.edu}

\date{\today}

\begin{abstract}
We use rank one Gaussian perturbations to derive a smooth stochastic approximation of the maximum eigenvalue function. We then combine this smoothing result with an optimal smooth stochastic optimization algorithm to produce an efficient method for solving maximum eigenvalue minimization problems, and detail a variant of this stochastic algorithm with monotonic line search. Overall, compared to classical smooth algorithms, this method runs a larger number of significantly cheaper iterations and, in certain precision/dimension regimes, its total complexity is lower than that of deterministic smoothing algorithms. 
\end{abstract}
\maketitle

% AA: no space before punctuation in English?

\section{Introduction}
We discuss applications of stochastic smoothing results to the design of efficient first-order methods for solving semidefinite programs. We focus here on the problem of minimizing the maximum eigenvalue of a matrix over a simple convex set $Q$ (the meaning of simple will be made precise later), i.e. we solve
\BEQ\label{eq:min-max-eig}
\min_{X\in Q} \lambdamax(X),
\EEQ
in the variable $X\in\symm_n$. Note that all primal semidefinite programs with fixed trace have a dual which can be written in this form \citep{Helm00b}. While moderately sized problem instances are solved very efficiently by interior point methods \citep{Bent01} with very high precision guarantees, these methods fail on most large-scale problems because the cost of running even one iteration becomes too high. When coarser precision targets are sufficient (e.g. spectral methods in statistical or geometric applications), much larger problems can be solved using first-order algorithms, which tradeoff a lower cost per iteration in exchange for a degraded dependence on the target precision.

So far, roughly two classes of first-order algorithms have been used to solve large-scale instances of the semidefinite program in~\eqref{eq:min-max-eig}. The first uses subgradient descent or a variant of the mirror-prox algorithm of \citep{Nemi79} that takes advantage of the geometry of $Q$ to minimize directly $\lambdamax(X)$. These methods do not exploit the particular structure of problem~\eqref{eq:min-max-eig} and need $\gO(D_Q^2/\epsilon^2)$ iterations to reach a target precision $\epsilon$, where $D_Q$ is the diameter of the set $Q$. Each iteration requires computing a leading eigenvector of the matrix $X$ at a cost of roughly $\gO(n^2 \log n)$ (see the Appendix for more details) and projecting $X$ on $Q$ at a cost written $p_Q$. Spectral bundle methods \citep{Helm00b} use more information on the spectrum of $X$ to speed up convergence, but their complexity is not well understood. More recently, \citep{Nest04a} showed that one could exploit the particular min-max structure of problem~\eqref{eq:min-max-eig} by first regularizing the objective using a ``soft-max'' exponential smoothing, then using optimal first-order methods for smooth convex minimization. These algorithms only require $\gO(D_Q\sqrt{\log n}/\epsilon)$ iterations, but each iteration forms a matrix exponential at a cost of $O(n^3)$. In other words, depending on problem size and precision targets, existing first-order algorithms offer a choice between two complexity bounds
\BEQ\label{eq:perf-bnds}
O\left(\frac{ D_Q^2 (n^2\log n+p_Q)}{\epsilon^2}\right) \quad \mbox{and} \quad O\left(\frac{D_Q\sqrt{\log n}(n^3+p_Q)}{\epsilon}\right).
\EEQ
Note that the constants in front of all these estimates can be quite large and actual numerical complexity depends heavily on the particular path taken by the algorithm, especially for adaptive variants of the methods detailed here (see \citep[\S 6]{Nest07} for an illustration on a simpler problem). In practice of course, these asymptotic worst case bounds are useful for providing general guidance in algorithmic choices, but remain relatively coarse predictors of performance for reasonable values of $n$ and $\epsilon$.

Many recent works have sought to move beyond these two basic complexity options. \citet{Over95} directly applied Newton's method to the maximum eigenvalue function, given a priori information on the multiplicity of this eigenvalue.  \citet{Bure03} and \citet{Jour08a} focus on instances where the solution is known to have low rank (e.g. matrix completion, combinatorial relaxations) and solve the problem directly over the set of low rank matrices. These formulations are nonconvex and their complexity cannot be explicitly bounded, but empirical performance is often very good. \citet{Lu07} focus on the case where the matrix has a natural structure (close to block diagonal). \citet{Judi08a} use a variational inequality formulation and randomized linear algebra to reduce the cost per iteration of first-order algorithms. Subsampling techniques were also used in \citep{dAsp11} to reduce the cost per iteration of stochastic averaging algorithms. Finally, in recent independent results similar in spirit to those presented here, \citet{Baes11} use stochastic approximations of the matrix exponential to reduce the cost per iteration of smooth first-order methods. The complexity tradeoff and algorithms in \citep{Baes11} are different from ours (roughly speaking, a $1/\epsilon$ term is substituted to the $\sqrt{n}$ term in our bound), but both methods seek to reduce the cost of smooth first-order algorithms for semidefinite programming using stochastic gradient oracles instead of deterministic ones. However, \citet{Baes11} use stochastic techniques to reduce the cost of computing classical smoothing steps (matrix exponential, etc.) and \citet{Judi08a} use them to reduce the cost of linear algebra operations. In this work, we directly use stochastic methods for smoothing.

In this paper, we use stochastic smoothing results, combined with an optimal accelerated algorithm for stochastic optimization recently developed by \citet{Lan09}, to derive a stochastic algorithm for solving~\eqref{eq:min-max-eig}. The algorithm detailed below requires only $O(\sqrt{n}/\epsilon)$ iterations, with each iteration computing a few sample leading eigenvectors of $(X+\epsilon \, \canonicalRandomVector\canonicalRandomVector^T/n)$ where $\canonicalRandomVector\sim{\mathcal N}(0,\idm_n)$. While in most applications of stochastic optimization the noise level is seen as exogenous, we use it here to control the tradeoff between number of iterations and cost per iteration. The algorithm requires fewer iterations than nonsmooth methods and has lower cost per iteration than smoothing techniques. In some configurations of the parameters~$(n,\epsilon,p_Q,D_Q)$, its total worst-case floating-point complexity is lower than that of both smooth and nonsmooth methods. Overall, the method has a cost per iteration comparable to that of nonsmooth methods while retaining some of the benefits of accelerated methods for smooth optimization.

The paper is organized as follows. In the next section, we briefly outline the stochastic smoothing algorithm for maximum eigenvalue minimization and compare its complexity with existing first-order algorithms. Section~\ref{s:smooth} details smoothing results on random rank one perturbations of the maximum eigenvalue function, highlighting in particular a phase transition in the spectral gap  depending on the spectrum of the original matrix. Section~\ref{s:sco} uses these smoothing results to produce a stochastic algorithm for maximum eigenvalue minimization, and describes an extension of the optimal stochastic optimization algorithm in \citep{Lan09} where the scale of the step size is allowed to vary adaptively (but monotonically). Section~\ref{s:speculative-math} informally discusses extensions of our results to other smoothing techniques, together with their impact on complexity. Section~\ref{s:numexp} presents some preliminary numerical experiments. An Appendix contains auxiliary material, including a detailed discussion of the cost of computing leading eigenpairs of a symmetric matrix, technical details about various functions that play a central role in our analysis, and a proof of the phase transition result for random rank-one perturbations.

\subsection*{Notation} Throughout the paper, we denote by $\lambda_i(X)$ the  eigenvalues of the matrix $X \in \symMatrices$, in decreasing order. For clarity, we will also use $\lambdamax(X)$ for the leading eigenvalue of $X$. When $\canonicalRandomVector$ denotes a vector in $\reals^n$, its $i$-th coordinate is denoted by $\coordCanonRandomVector_i$. We denote equality in law (for random variables) by $\stackrel{\mathcal{L}}{=}$ and $\weakCV$ stands for convergence in law. We use the notation $\gO_P$ with the standard probabilistic meaning (see \citep{vanderVaart98}, p.12). When we compute local Lipschitz constants, they are always computed with respect to Euclidian or Frobenius norms, unless otherwise noted. We call $\localLipConstant{\Gamma(X)}$ the local Lipschitz constant of the function $\Gamma$ at $X$.

\section{Stochastic smoothing algorithm}\label{s:summary}
We will solve a smooth approximation of problem~\eqref{eq:min-max-eig}, written
\BEQ\label{eq:smooth-min}
\BA{ll}
\mbox{minimize } & \smoothedF_k(X)\triangleq \Expect\left[\max_{i=1,\ldots,k} \lambdamax\left(X+\frac{\epsilon}{n}\canonicalRandomVector_i\canonicalRandomVector_{i}^T\right)\right]\\
\mbox{subject to} & X \in Q,\\
\EA\EEQ
in the variable $X\in\symm_n$, where  $Q\subset \symm_n$ is a compact convex set, $\canonicalRandomVector_i\iid {\mathcal N}(0, \idm_n )$, $\eps\geq 0$ is in $\reals$ and $k>0$ is a small constant (typically 3). We call $\smoothedF_k^*$ the optimal value of this problem. We also define $F_k(X)$ as the random valued function inside the expectation, with
\BEQ\label{eq:F_k}
F_k(X)\triangleq \max_{i=1,\ldots,k} \lambdamax\left(X+\frac{\epsilon}{n}\canonicalRandomVector_i\canonicalRandomVector_{i}^T\right)
\EEQ
so that $\smoothedF_k(X)= \Expect\left[F_k(X)\right]$. We have the following approximation result. 
\begin{lemma}\label{lemma:smoothingIsCkEpsUnifApproxToLambdaMax}
$\smoothedF_k(X)$ is a $c_k \epsilon$-uniform approximation of $\lambdamax(X)$, where 
$$
c_k = \Expect\left[\max_{i=1,\ldots,k} \|\canonicalRandomVector_i\|_2^2/n\right] \leq  \Expect\left[\textstyle \sum_{i=1}^k \|\canonicalRandomVector_i\|_2^2/n\right]=k\;.
$$
In other words, for all $X  \in \symm_n\,$ 
\begin{equation}\label{eq:smoothingIsCkEpsUnifApproxToLambdaMax}
\lambdamax(X)+\frac{\epsilon}{n}\leq \smoothedF_k(X)\leq  \lambdamax(X)+c_k \eps\;.
\end{equation}
\end{lemma}

\begin{proof}
We first establish the upper bound. The fact that $\lambdamax(\cdot)$ is subadditive on $\symMatrices$ gives 
$$
\lambdamax\left(X+\frac{\epsilon}{n}\canonicalRandomVector_i\canonicalRandomVector_{i}^T\right)\leq \lambdamax(X)+\lambdamax\left(\frac{\epsilon}{n}\canonicalRandomVector_i\canonicalRandomVector_{i}^T\right)=\lambdamax(X)+
\eps\frac{\norm{\canonicalRandomVector_i}^2}{n}\;,
$$
since $\lambdamax(\canonicalRandomVector_i\canonicalRandomVector_i\trsp)=\norm{\canonicalRandomVector_i}^2$.
%1-Lipschitz with respect to the spectral norm with $\lambdamax(\canonicalRandomVector_i\canonicalRandomVector_{i}^T)=\|\canonicalRandomVector_i\|_2^2$, yields
It follows that
\[
\max_{1\leq i \leq k} \lambdamax\left(X+\frac{\epsilon}{n}\canonicalRandomVector_i\canonicalRandomVector_{i}^T\right) \leq \max_{1\leq i \leq k} 
\lambdamax(X)+
\eps\max_{1\leq i \leq k}\frac{\norm{\canonicalRandomVector_i}^2}{n}\;, 
\]
and
\[
\smoothedF_k(X)=\Expect\left[\max_{i=1,\ldots,k} \lambdamax\left(X+\frac{\epsilon}{n}\canonicalRandomVector_i\canonicalRandomVector_{i}^T\right)\right] \leq \lambdamax(X) + c_k \epsilon\;.
\]
% where
% $$
% c_k = \Expect\left[\max_{i=1,\ldots,k} \|\canonicalRandomVector_i\|_2^2/n\right] \leq  \Expect\left[\textstyle \sum_{i=1}^k \|\canonicalRandomVector_i\|_2^2/n\right]=k \\
% $$
% depends only on $k$. 
Let us now prove the lower bound. The mapping $M\mapsto \lambdamax(X+M)$ is convex from $\symm_n$ to $\reals$ when $X\in \symm_n$. Therefore, Jensen's inequality applied to this mapping with the random varible $\canonicalRandomVector_i\canonicalRandomVector_i\trsp$ gives
$$
\lambdamax(X+\frac{\eps}{n}\Exp{\canonicalRandomVector_i\canonicalRandomVector_i\trsp})\leq \Exp{\lambdamax\left(X+\frac{\eps}{n}\canonicalRandomVector_i\canonicalRandomVector_i\trsp\right)}\;.
$$
Using $\Expect[\canonicalRandomVector_i\canonicalRandomVector_{i}^T]=\idm_n$, we conclude that 
\begin{align*}
\forall 1\leq i \leq k, \; \lambdamax(X+\frac{\epsilon}{n}\idm_n)&\leq \Expect\left[\lambdamax\left(X+\frac{\epsilon}{n}\canonicalRandomVector_i\canonicalRandomVector_{i}^T\right)\right]\;,  \text{ hence }\\
\lambdamax(X)+\frac{\epsilon}{n}&\leq \Expect\left[\max_{i=1,\ldots,k} \lambdamax\left(X+\frac{\epsilon}{n}\canonicalRandomVector_i\canonicalRandomVector_{i}^T\right)\right],
\end{align*}
which is the lower bound above.
\end{proof}

We begin by briefly introducing the smoothing results on~\eqref{eq:smooth-min} detailed in Section~\ref{s:smooth}, then describe our main algorithm.

\subsection{Smoothness of $\smoothedF_k(X)$} In Section~\ref{s:smooth}, we will show that the function $\smoothedF_k$ has a Lipschitz continuous gradient w.r.t. the Frobenius norm, i.e.
\[
\|\nabla \smoothedF_k(X) - \nabla \smoothedF_k(Y) \|_{F} \leq L \|X-Y\|_F\;,
\]
with (uniform) constant $L$ satisfying
\BEQ\label{eq:Lip}
L \leq C_k \frac{n}{\epsilon}\;,
\EEQ
where $C_k>0$ depends only on $k$ and is bounded whenever $k\geq 3$. We will see in Section~\ref{s:smooth} that this bound is quite conservative and that much better regularity is achieved when the spectrum of $X$ is well-behaved (see Theorem \ref{thm:PhaseTransition}).

\subsection{Gradient variance}
Section~\ref{s:smooth} also produces an explicit expression for the gradient of $\smoothedF_k$. Let $\phi_{i_0}$ be a leading eigenvector of the matrix $X+\frac{\epsilon}{n}\canonicalRandomVector_{i_0}\canonicalRandomVector_{i_0}^{T}$
where
$$
i_0=\argmax_{i=1,\ldots,k} \lambdamax\left(X+\frac{\epsilon}{n}\canonicalRandomVector_i\canonicalRandomVector_{i}^T\right).
$$
We will see that $i_0$ is unique with probability one and that we have
\BEQ\label{eq:var-bound}
\nabla \smoothedF_k(X)=\Expect\left[\phi_{i_0}\phi_{i_0}^T\right]
\quad \mbox{and} \quad
\Expect\left[\left\|\phi_{i_0}\phi_{i_0}^T-\nabla \smoothedF_k(X) \right\|_F^2\right] \leq 1\;.
\EEQ
Therefore the variance of the stochastic gradient oracle $\phi_{i_0}\phi_{i_0}^T$ is bounded by one. Once again, we will see in Section~\ref{s:smooth} that this bound too is often quite conservative.

\subsection{Stochastic algorithm}
Given an unbiased estimator for $\nabla \smoothedF_k$ with unit variance, the optimal algorithm for stochastic optimization derived in \citep{Lan09} will produce a (random) matrix $X_N$ such that 
\BEQ\label{eq:lan-N}
\Expect[\smoothedF_k(X_N)-\smoothedF_k^*]\leq \frac{4LD^2_Q}{\alpha N^2} + \frac{4 D_Q}{\sqrt{Nq}}
\EEQ
after $N$ iterations \citep[Corollary\,1]{Lan09}, where $L\leq C_k n/\epsilon$ is the Lipschitz constant of $\nabla \smoothedF_k$ discussed in the previous section, $\alpha$ is the strong convexity constant of the prox function, and $q$ is the number of independent sample matrices $\phi\phi^T$ averaged in approximating the gradient. Once again, we write $D_Q$ the diameter of the set $Q$ (see below for a precise definition) and $p_Q$, which appears in Table~\ref{tab:complex}, the cost of projecting a matrix $X\in\symm_n$ on the set $Q$.

Setting $N=2D_Q\sqrt{n}/\epsilon$ and $q=\lceil\max\{1,D_Q/(\epsilon \sqrt{n})\}\rceil$, the approximation bounds in Proposition~\ref{prop:smoothness} will then ensure $\Expect[\smoothedF_k(X_N)-\smoothedF_k^*]\leq 5\epsilon$. We compare in Table~\ref{tab:complex} the computational cost of the smooth stochastic algorithm in \citep[Corollary\,1]{Lan09} in this setting with that of the smoothing technique in \citep{Nest04a} and the nonsmooth stochastic averaging method. Recall that the cost of computing one leading eigenvector of $X+\canonicalDetVector\canonicalDetVector^T$ is of order $O(n^2\log n)$ (cf. Appendix) while that of forming the matrix exponential $\exp(X)$ is $O(n^3)$ \citep{Mole03}.

\begin{table}%\label{tab:complex}
\begin{center}
\extrarowheight 1.5ex
\begin{tabular}{r|c|c}
{\bf Algorithmic complexity} & Num. of Iterations & Cost per Iteration \\
\hline
Nonsmooth & $O\left(\frac{D_Q^2}{\epsilon^2}\right)$
  & $O(p_Q+n^2\log n)$ \\
Stochastic Smoothing & $O\left(\frac{D_Q\sqrt{n}}{\epsilon}\right)$ & $O\left(p_Q+\max\left\{1,\frac{D_Q}{\epsilon \sqrt{n}}\right\}n^2\log n\right)$\\
Deterministic Smoothing & $O\left(\frac{D_Q\sqrt{\log n}}{\epsilon}\right)$ & $O(p_Q+n^3)$\\
\end{tabular}
\caption{Worst-case computational cost of the smooth stochastic algorithm detailed here, the smoothing technique in \citep{Nest04a} and the nonsmooth subgradient descent method.
\label{tab:complex}}
\end{center}
\end{table}

Table~\ref{tab:complex} shows a clear tradeoff in this group of algorithms between the number of iterations and the cost of each iteration. In certain regimes for $(n,\epsilon)$, the total worst-case complexity of algorithm~\ref{alg:stoch} detailed on page~\pageref{alg:stoch} is lower than that of both smooth and nonsmooth methods. This is the case for instance when $D_Q\geq \sqrt{n}\eps$  and
\[
c_1 \max\left\{1,\frac{D_Q}{\epsilon \sqrt{n}}\right\} n^2 \log n \leq p_Q \leq c_2 n^{5/2} \sqrt{\log n}\;,
\]
for some absolute constants $c_1,c_2>0$. In practice of course, the constants in front of all these estimates can be quite large and the key contribution of the algorithm detailed here is to preserve some of the benefits of smooth accelerated methods (e.g. fewer iterations), while requiring a much lower computational (and memory) cost per iteration by exploiting the very specific structure of the $\lambdamax(X)$ function.

\section{Efficient Stochastic smoothing} \label{s:smooth}
In this section, we show how to regularize the function $\lambdamax(X)$ using stochastic smoothing arguments. We start by recalling a classical argument about Gaussian regularization and then improve smoothing performance by using explicit structural results on the spectrum of rank one updates of symmetric matrices.

\subsection{Gaussian smoothing}
The following is a standard result on Gaussian smoothing which does not exploit any structural information on the function $\lambdamax(X)$ except its Lipschitz continuity.
\begin{lemma}\label{lem:gauss-smooth}
Suppose $\mathsf{f}:\reals^m\rightarrow \reals$ is Lipschitz continuous with constant $\mu$ with respect to the Euclidean norm. The function $\mathsf{sf}$ such that
$$
\mathsf{sf}(x)=\Expect[\mathsf{f}(x+\mysigma \canonicalRandomVector)]\;,
$$
where $\canonicalRandomVector\sim\mathcal{N}(0,\idm_m)$ and $\mysigma>0$, has a Lipschitz continuous gradient with
\[
\|\nabla \mathsf{sf}(x) - \nabla \mathsf{sf}(y)\|_2\leq \frac{2\sqrt{m}\mu}{\mysigma}\|x-y\|_2.
\]
\end{lemma}
\begin{proof}
See \cite{Nest11} for a short proof and applications in gradient-free optimization.
\end{proof}

Let us consider the function $\smoothedF_{GUE}(X)$ taking values
\[
\smoothedF_{GUE}(X)=\Expect[\lambdamax(X+(\epsilon/\sqrt{n}) U)]\;,
\]
where $U\in\symm_n$ is a symmetric matrix with standard normal upper triangle coefficients.
Using convexity and positive homogeneity of the $\lambdamax(X)$ function, together with the fact that it is 1-Lipschitz with respect to the spectral norm and bounds on the largest eigenvalue of $U$ (which follow easily from either
\cite{Trotter84} or \cite{DavidsonSzarek01}), we see that this function is an $\epsilon$-approximation of $\lambdamax(X)$. 

Lemma \ref{lem:gauss-smooth} above shows that $\smoothedF_{GUE}(X)$ has a Lipschitz continuous gradient with constant bounded by $O\left({n^{3/2}}/{\epsilon}\right)$, since, with the notation of Lemma \ref{lem:gauss-smooth}, $m=n^2$. This approach was used e.g. in \citep{dAsp08c} to reduce the cost per iteration of a smooth optimization algorithm with approximate gradient, and by \citep{Nest11} to derive explicit complexity bounds on gradient free optimization methods. We present a short discussion on a finer bound on the Lipschitz-constant of this function in Section \ref{s:speculative-math}.

\subsection{Gradient smoothness} \label{ss:smoothness}
We recall the following classical result (which can be derived from results in \citep{Kato95} and \citep{LewisSendov02TwiceDiffSpectFunctions} and is proved in the Appendix for the sake of completeness) showing that the gradient of $\lambdamax(X)$ is smooth when the largest eigenvalue of $X$ has multiplicity one, with (local) Lipschitz constant controlled by the spectral gap.

\begin{theorem} \label{th:lip-gap}
Suppose $X\in\symm_n$ and call $\{\lambda_i(X)\}_{i=1}^n$ the decreasingly ordered eigenvalues of $X$. Suppose also that $\lambdamax(X)$, the largest eigenvalue of $X$, has multiplicity one. Let $Y$ be a symmetric matrix with $\norm{Y}_F=1$ and call
$$
\gamma(X,Y)=\lim_{t\tendsto 0} \frac{\partial^2 \lambdamax(X+tY)}{\partial t^2}.
$$
Call $\lambdamax$ the mapping $X\mapsto \lambdamax(X)$.
Then the local Lipschitz constant - with respect to the Frobenius norm - of the gradient of the mapping $\lambdamax$  is given by
\begin{equation} \label{eq:valueLipConstant}
\localLipConstant{\nabla \lambdamax(X)}=\sup_{Y\in \symm_n, \norm{Y}_F=1} \gamma(X,Y)=\frac{1}{\lambdamax(X)-\lambda_2(X)}.
\end{equation}
\end{theorem}
This result shows that to produce smooth approximations of the function $\lambdamax(X)$ using random perturbations, we need these perturbations to increase the spectral gap by a sufficient margin. We will see below that, up to a small trick, random rank one Gaussian perturbations of the matrix $X$ will suffice to achieve this goal.

\subsection{Rank one updates}\label{ss:rank-one}
The following proposition summarizes the information we will need about the impact of rank-one updates on the largest eigenvalue of a symmetric matrix. Equation \eqref{eq:gap-lb} below will prove useful later to control the smoothness of $\nabla \smoothedF_k(X)$.
\begin{proposition}\label{prop:impactLowRankUpdateOnTopEigenvalue}
Suppose $X\in\symm_n$ and has spectral decomposition $X=\orthoMat_X\trsp D_X \orthoMat_X$. Let $\canonicalDetVector\neq 0$ be a vector in $\reals^n$ which is not an eigenvector of $X$. Let $\mysigmasquared>0$ be in $\reals$. Then, $\lambdamax(X+(\mysigmasquared/n) \canonicalDetVector\canonicalDetVector\trsp)$ has multiplicity~1 and  $\lambdamax(X+(\mysigmasquared/n) \canonicalDetVector\canonicalDetVector\trsp)-\lambdamax(X)>0$. Let us call $\lambda_2$ the second largest eigenvalue of a symmetric matrix. Then, if $(\orthoMat_X\canonicalDetVector)_1$ is the first coordinate of the vector $\orthoMat_X\canonicalDetVector$, we have
\BEQ\label{eq:gap-lb}
\frac{\mysigmasquared (\orthoMat_X\canonicalDetVector)_1^2}{n} \leq \lambdamax\left(X+\frac{\mysigmasquared}{n} \canonicalDetVector\canonicalDetVector\trsp\right)-\lambdamax(X) \leq  \lambdamax\left(X+\frac{\mysigmasquared}{n} \canonicalDetVector\canonicalDetVector\trsp\right)-\lambda_2\left(X+\frac{\mysigmasquared}{n} \canonicalDetVector\canonicalDetVector^T\right)\;.
\EEQ
\end{proposition}
\begin{proof}
For $X\in\symm_n$, we call $\lambda(X)\in\reals^n$ the spectrum of the matrix $X$, in decreasing algebraic order. Whenever $\canonicalDetVector\neq 0$ is not an eigenvector of $X$ and $\mysigmasquared>0$, the leading eigenvalue $l_1$ of the matrix $X+(\mysigmasquared/n) \canonicalDetVector\canonicalDetVector\trsp$, is always strictly larger than $\lambda_1(X)$ \citep[see][\S8.5.3]{Golu90} and we write $l_1=\lambda_1(X)+t$, $t\geq 0$. Our aim is now to characterize $t$ and understand its properties. We note that 
$$
X+(\mysigmasquared/n) \canonicalDetVector\canonicalDetVector\trsp=\orthoMat_X\trsp \left[D_X+(\mysigmasquared/n) (\orthoMat_X\canonicalDetVector)(\orthoMat_X\canonicalDetVector)\trsp\right]\orthoMat_X\;.
$$
Since we are interested in eigenvalues, we assume without loss of generality that $X$ is diagonal. If $X$ were not diagonal, we would just need to replace $\canonicalDetVector$ by $(\orthoMat_X\canonicalDetVector)$ in what follows for all the statements to hold.

It is a standard result (see e.g Theorem 8.5.3 in \citep[see][\S8.5.3]{Golu90}) that the variable~$t$ is the unique positive solution of the {\em secular equation}
\BEQ\label{eq:secular}
s(t) \triangleq \frac{n}{\mysigmasquared}-\frac{\coordCanonDetVector_1^2}{t}-\sum_{i=2}^n\frac{\coordCanonDetVector_i^2}{(\lambda_1(X)-\lambda_i(X))+t}=0,
\EEQ
where $\coordCanonDetVector_i$ are the coefficients of the vector $v$; we give an elementary derivation of this result in Subsection~\ref{ss:secularEqn} in the Appendix. We plot the function $s(\cdot)$ for a sample matrix $X$ in Figure~\ref{fig:secular}. 
\begin{figure}[ht!]
\begin{center}
\psfrag{flambda}[b][t]{$s(t)$}
\psfrag{lambda}[t][b]{$\lambda_1(X)+t$}
\psfrag{lo}[c][b]{\small{$~~~\lambda_1(X)$}}
\psfrag{lp}[c][c]{\small{$\qquad\quad\lambda_1(X)+t^*$}}
\includegraphics[width=0.5 \textwidth]{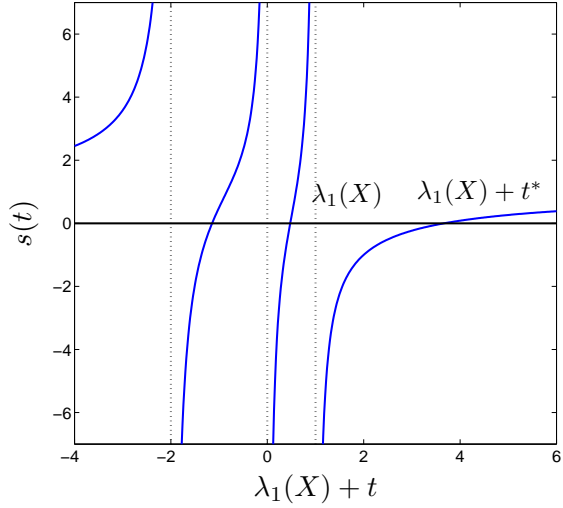}
\caption{Plot of $s(t)$ versus $\lambda_1(X)+t$. The matrix has dimension four and its spectrum is here $\{-2,-2,0,1\}$. The three leading eigenvalues of $X+\mysigmasquared \canonicalDetVector\canonicalDetVector^T$ are the roots of $s(t)$, the fourth eigenvalue is at -2.\label{fig:secular}}
\end{center}
\end{figure}

% Equation \eqref{eq:secular} implies that we have almost explicit expressions for the eigenvalue decomposition of the matrix
% $$
% X+\frac{\mysigmasquared}{n} \canonicalDetVector\canonicalDetVector\trsp
% $$
% where $\canonicalDetVector\in\reals^n$ and $\mysigma>0$. 

Having assumed that $X$ is diagonal, \citet[Th.\,8.5.3]{Golu90} also shows that if $\coordCanonDetVector_i\neq 0$ for $i=1,\ldots,n$ and $\mysigmasquared>0$, then $t>0$ and the eigenvalues of $X$ and $X+(\mysigmasquared/n) \canonicalDetVector\canonicalDetVector^T$ are interlaced, i.e.
\[
\lambda_n(X) \leq \lambda_n(X+\frac{\mysigmasquared}{n} \canonicalDetVector\canonicalDetVector^T)\leq \ldots \leq\lambda_2(X+\frac{\mysigmasquared}{n} \canonicalDetVector\canonicalDetVector^T) \leq \lambdamax(X) < \lambdamax(X+\frac{\mysigmasquared}{n} \canonicalDetVector\canonicalDetVector^T).
\]
This implies in particular that $\lambdamax(X+\frac{\mysigmasquared}{n} \canonicalDetVector\canonicalDetVector^T)$ has multiplicity 1. By construction, the function
\[
s^+(t)\triangleq\frac{n}{\mysigmasquared}-\frac{\coordCanonDetVector_1^2}{t}
\]
is an upper bound on $s(t)$ on the interval $(0,\infty)$. Since both functions are non-decreasing, the positive root of the equation $s^+(t)=0$ is a {lower bound} on the positive root $t^*$ of the equation $s(t)=0$. We therefore have 
\[
t^*\geq \frac{\mysigmasquared \coordCanonDetVector_1^2}{n}\;.
\]
Using interlacing, we also have
\[
\lambda_2(X+\frac{\mysigmasquared}{n} \canonicalDetVector\canonicalDetVector^T) \leq \lambda_1(X) \leq \lambda_1(X)+t^* = \lambda_1(X+\frac{\mysigmasquared}{n} \canonicalDetVector\canonicalDetVector^T).
\]
This gives a lower bound on the spectral gap of the perturbed matrix
$$
\frac{\mysigmasquared \coordCanonDetVector_1^2}{n} \leq t^* \leq  \lambda_1(X+\frac{\mysigmasquared}{n} \canonicalDetVector\canonicalDetVector^T)-\lambda_2(X+\frac{\mysigmasquared}{n} \canonicalDetVector\canonicalDetVector^T)\;,
$$
which yields~\eqref{eq:gap-lb} and will allow us to control the smoothness of $\nabla \smoothedF_k(X)$.
\end{proof}

\subsection{Low rank Gaussian smoothing}\label{ss:gauss-rk-one}
We now come back to the objective function of Problem~\eqref{eq:smooth-min}, written
\[
\smoothedF_k(X)= \Expect\left[\max_{i=1,\ldots,k} \lambdamax\left(X+\frac{\mysigmasquared}{n}\canonicalRandomVector_i\canonicalRandomVector_{i}^T\right)\right]\;,
\]
where $\canonicalRandomVector_i$ are i.i.d. ${\mathcal N}(0, \idm_n )$ and $k>0$ is a small constant. We first show that we can differentiate under the expectation in the definition of $\smoothedF_k(X)$. This requires a few preliminaries which we now present.

\begin{lemma}\label{lem:density}
Let $\lambda_1(X)+T$ be the largest eigenvalue of the matrix $X+(\epsilon/n) \canonicalRandomVector\canonicalRandomVector^T$, where $X \in \symm_n$ is a given deterministic matrix and $\canonicalRandomVector\sim {\mathcal N}(0,\idm_n)$. Then the random variable $T$ has a density on $[0,\infty)$.
\end{lemma}
The proof of this lemma is in the Appendix in~\S\ref{subsubsec:densityForT}. 
Two corollaries immediately follow. The first one shows that two perturbed eigenvalues obtained from independent rank one perturbations are different with probability one.

\begin{corollary}\label{coro:eigsAreDiffWithProba1WhenGaussianSmoothing}
Suppose $l_{1,1}=\lambdamax(X+(\epsilon/n) \canonicalRandomVector_1\canonicalRandomVector_1^{T})$ and $l_{1,2}=\lambdamax(X+(\epsilon/n) \canonicalRandomVector_2\canonicalRandomVector_2^{T})$, where $\canonicalRandomVector_1$ and $\canonicalRandomVector_2$ are independent with distribution ${\mathcal N}(0,\idm_n)$.
Then $l_{1,1}\neq l_{1,2}$ with probability one.
\end{corollary}
\begin{proof}
The result follows from Lemma~\ref{lem:density} since $l_{1,1}-\lambdamax(X)$ and $l_{1,2}-\lambdamax(X)$ are two independent draws from a distribution with a density on $[0,\infty)$ and $P(l_{1,1}-\lambdamax(X)=0)=P(l_{1,2}-\lambdamax(X)=0)=0$.
\end{proof}

The second corollary shows that the maximum of (independent) perturbed eigenvalues is differentiable with probability one and bounds its Lipschitz constant.

\begin{corollary}\label{cor:diff}
Let $X\in \symm_n$ and suppose $l_{1,i}=\lambdamax(X+(\epsilon/n) \canonicalRandomVector_i\canonicalRandomVector_i^{T})$, where $\canonicalRandomVector_i \iid {\mathcal{N}}(0,\idm_n)$ for $i=1,\ldots,k$. The mapping $F_k: X\rightarrow \max_{i=1,\ldots,k}l_{1,i}$ is differentiable with probability one. Then, if $i_0=\argmax_{1\leq i \leq k} l_{1,i}$ and $\phi_{i_0}$  is an eigenvector associated with the eigenvalue $l_{1,i_0}$, its gradient is
\begin{equation}\label{eq:gradientFk}
\nabla F_k(X)=\phi_{i_0}\phi_{i_0}\trsp\;.
\end{equation}
Also, with probability 1, the local Lipschitz constant of $\nabla F_k$ is bounded by 
\begin{equation}\label{eq:boundLipConstantFkDeterministicVersion}
\localLipConstant{\nabla F_k(X)}\leq \frac{1}{F_k(X)-\lambdamax(X)}\;.
%\localLipConstant{F_k(X)}\leq \frac{1}{2}\frac{1}{\lambdamax(X+(\epsilon/n) \canonicalRandomVector_{i_0}\canonicalRandomVector_{i_0}^{T})-\lambdamax(X)}\;,
\end{equation}

\end{corollary}
\begin{proof}
We first recall that it is well-known (and indeed follows from results in \citep{Kato95}) that if a matrix $M_0$ has a unique largest eigenvalue, the gradient of $M\mapsto \lambdamax(M)$ at $M_0$ is simply $\phi_0\phi_0\trsp$, where $\phi_0$ is an eigenvector associated with $\lambdamax(M_0)$. 

Corollary \ref{coro:eigsAreDiffWithProba1WhenGaussianSmoothing} shows that with probability 1, there exists a unique $i_0$ such that $l_{1,i_0}=F_k(X)$. Furthermore, since with probability 1, $\canonicalRandomVector_{i_0}$ is not an eigenvector of $X$, Proposition \ref{prop:impactLowRankUpdateOnTopEigenvalue} shows that the multiplicity of the largest eigenvalue of $X+(\epsilon/n) \canonicalRandomVector_{i_0}\canonicalRandomVector_{i_0}^{T}$ is one. This implies that $X\mapsto \lambdamax(X+(\epsilon/n) \canonicalRandomVector_{i_0}\canonicalRandomVector_{i_0}^{T})$ is differentiable at $X$ with probability 1. 
Lemma \ref{lemma:diffSupSeveralFunctions} then applies and shows that $F_k$ is differentiable at $X$ with probability~1. Our reminder on the gradient of $M\mapsto \lambdamax(M)$ gives the value of the differential. The last part of the corollary follows from Lemma \ref{lemma:localLipConstantLargestEigSeveralMatrices}, whose assumptions are clearly satisfied with probability~1.
\end{proof}

We now use these preliminary results to prove the main result of this section, namely a bound on the Lipschitz constant of the gradient of $\smoothedF_k(X)$ defined above, using the spectral gap bound in~\eqref{eq:gap-lb}.

\begin{proposition}\label{prop:smoothness}
Let $\left\{\canonicalRandomVector_{i}\right\}_{i=1}^k$ be i.i.d. ${\mathcal N}(0,\idm_n)$, $k\geq 3$ be an integer and $X \in \symm_n$. The function $\smoothedF_k$ such that
\[
\smoothedF_k(X)=\Expect\left[\max_{i=1,\ldots,k} \lambdamax(X+(\epsilon/n)\canonicalRandomVector_i\canonicalRandomVector_i^{T})\right]
\]
is smooth. The Lipschitz constant $L$ of its gradient w.r.t. the Frobenius norm satisfies
\[
L\leq C_k \frac{n}{\epsilon}
\quad \mbox{where} \quad
C_k=\frac{k}{k-2}\;.
\]
%when $k\geq 3$.
\end{proposition}
\begin{proof}
The fact that $\smoothedF_k$ is smooth follows from Equation \eqref{eq:gradientFk} and the fact that we can interchange expectation and differentiation here. Details about the validity of this interchange - whose proof requires care -  are in the Appendix in Lemma \ref{lemma:interchangeDiffAndEexpect}. We assume without loss of generality that $X$ is diagonal - Lemma~\ref{lemma:impactOfRotaInvariance} proving that we can do so. Let us call $\coordCanonRandomVector_{i}$ the first coordinate of the vector $\canonicalRandomVector_{i}$. The spectral gap bound in~\eqref{eq:gap-lb} gives 
$$
\forall \, i, 1\leq i \leq k\;, \; \; \;\lambdamax(X+\frac{\eps}{n}\canonicalRandomVector_{i} \canonicalRandomVector_{i}\trsp)-\lambdamax(X)\geq \frac{\eps}{n}\coordCanonRandomVector_{i}^2\;.
$$
It follows that 
\begin{align*}
F_k(X)-\lambdamax(X)&\geq \frac{\eps}{n} \max_{1\leq i \leq k} \coordCanonRandomVector_{i}^2\;, \text{ and }\\
\frac{1}{F_k(X)-\lambdamax(X)}&\leq \frac{n}{\epsilon} \frac{1}{\max_{i=1,\ldots,k} \coordCanonRandomVector_{i}^2}\;.
% \frac{n}{\eps}\min_{1\leq i \leq k} \frac{1}{\coordCanonRandomVector_{i}^2}\;.
\end{align*}
The results of Corollary \ref{cor:diff} then guarantee that with probability 1, we have
$$
\localLipConstant{\nabla F_k(X)} \leq \frac{n}{\epsilon} \frac{1}{\max_{i=1,\ldots,k} \coordCanonRandomVector_{i}^2}\;,
% =\frac{n}{\epsilon} \min_{i=1,\ldots,k}~\frac{1}{\coordCanonRandomVector_{i}^2}\;.
$$
and therefore
\[
\localLipConstant{\nabla \smoothedF_k(X)} \leq \Expect \left[ \frac{n}{\epsilon} \frac{1}{\max_{i=1,\ldots,k} \coordCanonRandomVector_{i}^2} \right] \leq 
\Expect \left[ \frac{n}{\epsilon}\,\frac{1}{\sum_{i=1}^k \coordCanonRandomVector_{i}^2/k} \right]
=\Expect \left[ \frac{n}{\epsilon}\,\frac{k}{\chi_k^2} \right]
\]
where $\chi_k^2$ is $\chi^2$ distributed with $k$ degrees of freedom. The fact that
\[
\Expect\left[\frac{1}{\chi_k^2}\right]=\frac{1}{k-2}% \frac{1}{2^{k/2}\Gamma(k/2)}\int_0^\infty t^{\frac{k-2}{2}-1}e^{-t/2}dt=\frac{\Gamma\left(\frac{k}{2}-1\right)}{\Gamma\left(\frac{k}{2}\right)2},
\]
whenever $k\geq 3$ - see e.g. \citep[p. 487]{MardiaKentBibby79} - yields
$$
\forall X \in \symm_n\;, \; \; \localLipConstant{\nabla \smoothedF_k(X)} \leq C_k \frac{n}{\mysigmasquared}\;.
$$
The function $\nabla \smoothedF_k$ is thus Lipschitz with Lipschitz constant 
$$
L\leq \sup_{X \in \symm_n} \localLipConstant{\nabla \smoothedF_k(X)} \leq C_k \frac{n}{\mysigmasquared}\;,
$$
which concludes the proof.
\end{proof}

\noindent Note that the bound above is a bit coarse; numerical simulations show that for independent ${\mathcal N}(0,1)$ random variables $\{\coordCanonRandomVector_i\}_{i=1}^3$, $$\Expect\left[{1}/{\max\{\coordCanonRandomVector_1^2,\coordCanonRandomVector_2^2,\coordCanonRandomVector_3^2\}}\right]=1.5...$$ while $C_3=3$, for example. We could of course use the density of the minimum above to get a more accurate bound, but then $C_k$ would not have a simple closed form.

\subsection{Gradient variance}\label{ss:grad-var}
In this section, we will bound the variance of $\nabla F_k$, the stochastic gradient oracle approximating $\nabla \smoothedF_k$. 

% \bluetext{
% Let us call
% $$
% g(X,\canonicalRandomVector)=\lambdamax(X+\frac{\mysigmasquared}{n} \canonicalRandomVector\canonicalRandomVector\trsp).
% $$
% Because of the rotational invariance of both $\lambdamax(\cdot)$ and of the Gaussian distribution, we can assume without loss of generality that $X$ is diagonal and that its largest eigenvalue has multiplicity~$l$.
% }
\begin{lemma}\label{eq:grad-var}	
Let $X\in\symm_n$  and $\canonicalRandomVector_i\iid {\mathcal N}(0,\idm_n)$, the gradient of $\smoothedF_k(X)$
% \[
% \smoothedF_k(X)=\Expect\left[\max_{i=1,\ldots,k} \lambdamax(X+(\epsilon/n)\canonicalRandomVector_i\canonicalRandomVector_i^{T})\right]
% \]
 is given by
\begin{equation} \label{eq:gradientSmoothedFk}
\nabla \smoothedF_k(X) = \Expect[\phi_{i_0}\phi_{i_0}^{T}]
\end{equation}
where $\phi_{i_0}$ is the leading eigenvector of the matrix $X+\frac{\epsilon}{n}\canonicalRandomVector_{i_0}\canonicalRandomVector_{i_0}^{T}$, and
\[
i_0=\argmax_{i=1,\ldots,k} \lambdamax\left(X+\frac{\epsilon}{n}\canonicalRandomVector_i\canonicalRandomVector_i^{T}\right).
\]
We have
\BEQ\label{eq:var-bnd}
\Expect\left[\|\phi_{i_0}\phi_{i_0}^{T}-\Expect[\phi_{i_0}\phi_{i_0}^{T}]\|_F^2\right]= 1-\Tr\left(\nabla \smoothedF_k(X)^2\right) \leq 1,
\EEQ
where $\Tr\left(\nabla \smoothedF_k(X)\right)=1$ by construction.
\end{lemma}
\begin{proof}
Equation \eqref{eq:gradientSmoothedFk} follows from Equation \eqref{eq:gradientFk} and the fact that we can interchange expectation and differentiation here (see Lemma \ref{lemma:interchangeDiffAndEexpect} for details). We now focus on the variance
$
\Exp{\norm{\phi_{i_0} \phi_{i_0}\trsp-\Exp{\phi_{i_0} \phi_{i_0}\trsp}}_F^2}.
$
Recall that for any symmetric matrix $M$, $\norm{M}_F^2=\trace(M\trsp M)=\trace{M^2}$. 
The matrix $\phi_{i_0} \phi_{i_0}\trsp-\Exp{\phi_{i_0} \phi_{i_0}\trsp}$ is symmetric. So we can rewrite 
$$
\norm{\phi_{i_0} \phi_{i_0}\trsp-\Exp{\phi_{i_0} \phi_{i_0}\trsp}}_F^2=\trace{(\phi_{i_0} \phi_{i_0}\trsp-\Exp{\phi_{i_0} \phi_{i_0}\trsp})^2}\;.
$$
Using the fact that $\phi_{i_0}\trsp \phi_{i_0}=1$, we see that $(\phi_{i_0} \phi_{i_0}\trsp)^2=\phi_{i_0} \phi_{i_0}\trsp$. Therefore, 
$$
\Exp{\trace{(\phi_{i_0} \phi_{i_0}\trsp)^2}}=\Exp{\trace{\phi_{i_0}\phi_{i_0}\trsp}}=\Exp{\phi_{i_0}\trsp\phi_{i_0}}=1\;.
$$
Recalling that  $\Exp{\phi_{i_0} \phi\trsp_{i_0}}=\nabla \smoothedF_k(X)$, we have shown that $\trace{\left(\nabla \smoothedF_k(X)\right)}=1$. We also see that 
$$
\Exp{\trace{\left(\phi_{i_0} \phi_{i_0}\trsp-\Exp{\phi_{i_0} \phi_{i_0}\trsp}\right)^2}}=\Tr{\left(\nabla \smoothedF_k(X)\right)}-\trace{\left(\nabla \smoothedF_k(X)^2\right)}=1-\trace{\left(\nabla \smoothedF_k(X)^2\right)}\leq 1\;.
$$
which is the desired result.
\end{proof}

Furthermore, we show in Lemma~\ref{lem:sim-diag} in the Appendix that $\nabla \smoothedF_k$ is diagonalizable in the same basis as $X$. In particular, when $X$ is diagonal, so is $\nabla \smoothedF_k$. Simply using the fact that $\phi_{i_0}$ is an eigenvector, we have of course
\BEQ\label{eq:var-abs-bnd}
\|\phi_{i_0}\phi_{i_0}^{T}-\Expect[\phi_{i_0}\phi_{i_0}^{T}]\|_F^2\leq 4
\EEQ
which means that the gradient will naturally satisfy the ``light-tail'' condition A2 in \citep{Lan09} for $\sigma^2=4$. The bound in~\eqref{eq:var-bnd} together with the proof above (in particular Equation \eqref{eq:coordsEigenVecRankOnePerturb}) show that when the spectral gaps $\lambda_1(X)-\lambda_i(X)$ are large, the diagonal of $\nabla \smoothedF_k(X)$ is approximately sparse. In that scenario, $\Tr(\nabla \smoothedF_k(X)^2)$ is close to $\Tr(\nabla \smoothedF_k(X))$, hence close to one, and the variance of the gradient oracle is small.

\subsection{A phase transition} \label{ss:phase-trans}
We can push our analysis if the impact of the low rank perturbation a little bit further. We focus again on the properties of a random rank one perturbation of a deterministic matrix $X$, specifically $X(\mysigmasquared)=X+(\mysigmasquared/n)\canonicalRandomVector\canonicalRandomVector\trsp$, where $\canonicalRandomVector \sim {\mathcal N}(0,\idm_n)$. As we will see, the bounds we obtained above are quite conservative and the Lipschitz constant of the gradient is in fact much lower than $n/\epsilon$ when the spectrum of $X$ is well-behaved (in a sense that will be made clear later). In particular, we will observe that there is a {\em phase transition phenomenon} in $\mysigmasquared$. Let us call $T=\lambdamax(X(\mysigmasquared))-\lambdamax(X)$. If the perturbation scale $\mysigmasquared$ is small, $T$ is of order $1/n$ (the worst-case bound we obtained above). If $\mysigmasquared$ is large, $T$ is of order one. And if $\mysigmasquared$ has a critical value, then $T$ is $\gO_P(1/\sqrt{n})$.

The next theorem is asymptotic in nature but is informative in practice even for moderate size matrices. We make the dependence on $n$, the dimension of the matrices we are working with, explicit everywhere. This undoubtedly makes for somewhat cumbersome notations but also makes the statement of the results less ambiguous. We will work under the following assumptions:
\begin{itemize}
\item[A1] $X_n \in \symm_n$. Its eigenvalues are denoted $\lambda_1(n)\geq \lambda_2(n)\geq \ldots\geq \lambda_n(n)$. $\lambda_1(n)$ has multiplicity $l_n \in \mathbb{N}$. There exists a constant $l \in \mathbb{N}$ such that $l_n\leq l$ for all $n$. We call $\gamma_n=\lambda_1(n)-\lambda_{l_n+1}(n)$ and assume that there exists a constant $\gamma$ such that $\gamma_n\geq \gamma>0$. We call $\lambda_1(n)-\lambda_i(n)=\gamma_n+\delta_{i,n}$, for $i>l(n)$. Of course, $\delta_{i,n}\geq 0$. 
\item[A2] $\mysigmasquared_n$ is a sequence in $\mathbb{R}$. We assume that $\mysigmasquared_n\asymp 1$, i.e $\liminf_{n\tendsto \infty} \mysigmasquared_n >0$ and 
$\limsup_{n\tendsto \infty} \mysigmasquared_n <\infty$. 
\item[A3] We assume that there exists a constant $C$, independent of $n$ such that 
$$
\frac{1}{\gamma^2}> \frac{1}{n}\sum_{j=l_n+1}^n \frac{1}{(\gamma_n+\delta_{j,n})^2}>  \frac{1}{n}\sum_{j=l_n+1}^n \frac{1}{(\mysigmasquared_n+\gamma_n+\delta_{j,n})^2}>  C\;.
$$
\end{itemize}

\begin{theorem}[\textbf{Phase transition for the largest eigenvalue: rank one perturbation}]\label{thm:PhaseTransition}
Assume that Assumptions A1-A3 above are satisfied and consider the matrix
$$
X_n(\mysigmasquared_n)=X_n+\frac{\mysigmasquared_n}{n} \canonicalRandomVector\canonicalRandomVector\trsp\;, 
\text{ where }\canonicalRandomVector\sim \mathcal{N}(0,\id_n) \;. 
$$
Define $\mysigmasquared_{0,n}$ by
$$
\frac{1}{\mysigmasquared_{0,n}}=\frac{1}{n}\sum_{j=l_n+1}^n\frac{1}{\gamma_n+\delta_{j,n}}\;.
$$
Call, for i.i.d ${\mathcal N}(0,1)$ random variables $\{\coordCanonRandomVector_{j,n}\}_{j=1}^n$, $\chi^2_{l_n}=\sum_{j=1}^{l_n} \coordCanonRandomVector_{j,n}^2$,
$$
\xi_{1,n}=\frac{1}{\sqrt{n}}\sum_{j=l_n+1}^n \frac{\coordCanonRandomVector_{j,n}^2-1}{\gamma_n+\delta_{j,n}}=\gO_P(1)\; \text{ and }
\zeta_{1,n}=\frac{1}{n}\sum_{j=l_n+1}^n \frac{\coordCanonRandomVector_{j,n}^2}{(\gamma_n+\delta_{j,n})^2}=\gO_P(1)\;.
$$
We have the following three situations:
\begin{enumerate}
\item If $0<\mysigmasquared_n<\mysigmasquared_{0,n}$ and $\liminf_{n\tendsto \infty} [\mysigmasquared_{0,n}-\mysigmasquared_n]>0$, as $n\tendsto \infty$,
$$
\lambdamax[X_n(\mysigmasquared_n)]=\lambdamax[X_n]+\frac{W_{1,n}}{n}+\frac{W_{2,n}}{n^{3/2}}+\gO_P\left(\frac{1}{n^2}\right),
$$
where
\[
W_{1,n}=\frac{\chi^2_{l_n}}{1/\mysigmasquared_n-1/\mysigmasquared_{0,n}}
\quad \mbox{and} \quad
W_{2,n}=\frac{W_{1,n} \xi_{1,n}}{1/\mysigmasquared_n-1/\mysigmasquared_{0,n}}\;.
\]
\item If $\mysigmasquared_n=\mysigmasquared_{0,n}$, as $n\tendsto \infty$,
$$
\lambdamax[X_n(\mysigmasquared_n)]=\lambdamax[X_n]+\frac{W_{1,n}}{\sqrt{n}}+\gO_P\left(\frac{1}{n}\right),
$$
where
$$
W_{1,n}=\frac{\xi_{1,n}+\sqrt{\xi_{1,n}^2+4\chi^2_{l_n}\zeta_{1,n}}}{2\zeta_{1,n}}.
$$
\item If $\mysigmasquared_n>\mysigmasquared_{0,n}$ and $\liminf_{n\tendsto\infty} [\mysigmasquared_n-\mysigmasquared_{0,n}]>0$, call $t_{0,n}>0$, the (unique) positive solution of
$$
\frac{1}{\mysigmasquared_n}=\frac{1}{n}\sum_{j=l_n+1}^n\frac{1}{t_{0,n}+\gamma_n+\delta_{j,n}}.
$$
Note that $t_{0,n}\leq (1-l_n/n)\mysigmasquared_n$. Then, as $n\tendsto \infty$, 
$$
\lambdamax[X_n(\mysigmasquared_n)]=\lambdamax[X_n]+t_{0,n}+\frac{W_{1,n}}{\sqrt{n}}+\gO_P\left(\frac{1}{n}\right).
$$
Here, $W_{1,n}=\frac{\xi(t_{0,n})}{\zeta(t_{0,n})}\;,$ where
\begin{align*}
\xi(t_{0,n})&=\frac{1}{\sqrt{n}} \sum_{j=l_n+1}^n \frac{\coordCanonRandomVector_{j,n}^2-1}{t_{0,n}+\gamma_n+\delta_{j,n}}=\gO_P(1)\;, \text{ and }\\
\zeta(t_{0,n})&=\frac{1}{n} \sum_{j=l_n+1}^n \frac{1}{(t_{0,n}+\gamma_n+\delta_{j,n})^2}=\gO(1)\;.
\end{align*}
\end{enumerate}
\end{theorem}

\begin{proof} The strategy is the following. We are looking for the zeros of a certain random function - defined in the secular equation -  which can be seen as a perturbation of a deterministic function. Hence, it is natural to use ideas from asymptotic root finding problems \citep[see][pp. 36-43]{MillerAppliedAsymptoticAnalysis06}, to expand the solution in powers of the size of the perturbation. We note that a similar idea was used in \citep{NadlerAOS08}, which  focused on a different random matrix problem. We now turn to the proof.

\subsubsection{Preliminaries}\label{subsubsec:mainArg}
Let us call
\begin{align*}
g_{l_n,n}(t)&=\frac{1}{n}\sum_{j=l_n+1}^n\frac{1}{t+\gamma_n+\delta_{j,n}}\;,\\
h_{l_n,n}(t)&=\frac{1}{n}\sum_{j=l_n+1}^n\frac{\coordCanonRandomVector_{j,n}^2}{t+\gamma_n+\delta_{j,n}}\;,\\
% g(t)&=\frac{\sum_{j=1}^l \coordCanonRandomVector_j^2}{n}\frac{1}{t}+\frac{1}{n}\sum_{j=l_n+1}^n\frac{1}{t+\gamma+\delta_j}=\frac{\sum_{j=1}^l \coordCanonRandomVector_j^2}{n}\frac{1}{t}+g_l(t)\;,\\
h_n(t)&=\frac{\sum_{j=1}^{l_n} \coordCanonRandomVector_{j,n}^2}{n}\frac{1}{t}+h_{l_n,n}(t)\;.
\end{align*}
Recall that if $\lambdamax[X_n(\eps_n)]=\lambdamax[X_n]+T$, $T$ is the unique positive solution of the equation
\begin{equation}\label{eq:KeyEq}
\frac{1}{\mysigmasquared_n}=h_n(T)=\frac{\sum_{j=1}^{l_n} \coordCanonRandomVector_{j,n}^2}{n}\frac{1}{T}+\frac{1}{n}\sum_{j=l_n+1}^n\frac{\coordCanonRandomVector_{j,n}^2}{T+\gamma_n+\delta_{j,n}}\;.
\end{equation}
It is clear that $T\geq (\mysigmasquared_n/n)\sum_{j=1}^{l_n} \coordCanonRandomVector_{j,n}^2$. Also, $h_n'(t)<0$ on $(0,\infty)$, so $h_n$ is invertible. We note that
$$
\var\left(\frac{1}{n}\sum_{j=l_n+1}^n\frac{\coordCanonRandomVector_{j,n}^2-1}{t+\gamma_n+\delta_{j,n}}\right)=\frac{1}{n} \left[\frac{1}{n}\sum_{j=l_n+1}^n\frac{2}{(t+\gamma_n+\delta_{j,n})^2}\right]\leq \frac{2}{n}\frac{1}{\gamma^2}=\gO\left(\frac{1}{n}\right)\;.
$$
It therefore follows from Chebyshev's inequality that the error made when replacing $h_{l_n,n}$ by $g_{l_n,n}$ when seeking the root of Equation \eqref{eq:KeyEq} is $\gO_P(1/\sqrt{n})$.

Our strategy is to expand $T$ in powers (possibly non-integer) of $1/n$. If we can find an approximate solution $t(m)$ of Equation \eqref{eq:KeyEq}, such that
$$
|h_n(t(m))-\frac{1}{\mysigmasquared_n}|=\gO_P(n^{-\beta})\;, \text{ for some } \beta \;,
$$
we claim that
$$
|t(m)-T|=\gO_P(n^{-\beta})\;.
$$
This is because $h_n$ is, at $\coordCanonRandomVector_{j,n}$ fixed, a Lipschitz function on $(\frac{\mysigmasquared_n\sum_{j=1}^{l_n} \coordCanonRandomVector_{j,n}^2}{n},\infty)$, and its Lipschitz constant is bounded below with high-probability on any compact subinterval of this interval. Hence,
we have, if $\norm{h_n^{-1}}_{\textrm{L},t(m),T}$ is the Lipschitz constant of $h_n^{-1}$ over an interval to which both $t(m)$ and $T$ belong, 
$$
|t(m)-T|=|h_n^{-1}(h_n(t(m)))-h_n^{-1}(h_n(T))|\leq \norm{h_n^{-1}}_{\textrm{L},t(m),T} |h_n(t(m))-\frac{1}{\mysigmasquared_n}|=\gO_P(n^{-\beta})\;.
$$

Note that if we can show that $|h_n'(y)|>C n^{b}$ in an interval containing both $t(m)$ and $T$, then $\norm{h_n^{-1}}_{\textrm{L},t(m),T}\leq n^{-b}/C$ and we  get by the same token
$$
|h_n(t(m))-\frac{1}{\mysigmasquared_n}|=\gO_P(n^{-\beta}) \Longrightarrow |t(m)-T|=\gO_P(n^{-(\beta+b)})\;.
$$
More details about these estimates are given in \ref{subsubsec:detailedTreatmentPhaseTFirstCase}, where we carry out a detailed proof.

To summarize, if we can come up with $t(m)$ which is a near solution of the equation $h_n(t)=1/\mysigmasquared_n$, it will be a good approximation of $T$. The quality of the approximation is detailed in the estimates above. In the proof below, we will exhibit such $t(m)$'s and, from them, get fine approximations of $T$. That is our strategy.

We finally recall that by definition
$$
\frac{1}{\mysigmasquared_{0,n}}=g_{l_n,n}(0)=\frac{1}{n}\sum_{j=l_n+1}^n\frac{1}{\gamma_n+\delta_{j,n}}\;.
$$

\textbf{Intuitive explanations} The following might help in giving the reader a sense of where the results come from. We wish to find an approximation of the root $T$ of the equation 
$1/\eps_n=h_n(T)$. The overall strategy is to write a Laurent-series expansion of $h_n$ around $x_n$, a real such that $h_n(x_n)-1/\eps_n$ is small. Practically, calling $A_{k,x_n}(h_n)$ our expansion to order $k$ of $h_n$ around $x_n$, we solve exactly the equation $A_{k,x_n}(h_n)(t)=1/\eps_n$. This strategy amounts practically to dropping various $\gO_P$ terms in our expansions of $h_n$ and solving the corresponding equations. Let us call $x^*_n$ the solution of $A_{k,x_n}(h_n)(t)=1/\eps_n$. Our proof shows that $T$ is indeed close to $x^*_n$, to various order of accuracy.

More specifically, we break $h_n$ into a component that stays bounded when $t\tendsto 0$ - this is what $h_{l_n,n}$ is - and a component that behaves like $1/(nt)$ as $t\tendsto 0$. 

\textbf{Cases 1) and 2) of the Theorem} In these cases, it is clear that if $c>0$, $\lim_{t\tendsto c} h_{n}(t)<1/\eps_n$ with high-probability. This suggests that $T\tendsto 0$ with high-probability. Hence our strategy is to write $h_{l_n,n}(t)=P_{l_n,n}(t)+\gO_P(t^{\alpha})$, where $P_{l_n,n}$ is a polynomial and $\alpha$ an integer, i.e expand $h_{l_n,n}$ in powers of $t$ for $t$ close to 0, and instead of solving $h_n(T)=1/\eps_n$, solve the approximating equation $h_n(x)-h_{l_n,n}(x)+P_{l_n,n}(x)=1/\eps_n$. This simply amounts to dropping the $O_P(t^{\alpha})$ term from our (Laurent-series) expansion of $h_n(t)$ in a neighborhood of $0$. This latter equation is a polynomial equation - hence it is easy to solve. Call $x^*_n$ its solution. By construction, it is fairly clear that $x^*_n$ is such that $h_n(x^*_n)$ is close to $1/\eps_n$. The proof makes this statement fully rigorous and pushes further to give rigorous statements concerning $T-x_n^*$, which is really the quantity we are interested in. 

\textbf{Case 3) of the Theorem} In this case, it is clear that $T$ has to remain bounded away from 0, since $h_{l_n,n}(0)>1/\eps_n$ with probability going to 1. Hence, we employ the same strategy as the one described above, except that we expand $h_n(t)$ around $t_{0,n}$, a non-random sequence bounded away from 0 picked such that $h_n(t_{0,n})-1/\eps_n\tendsto 0$ in probability. $h_n$ is linearized around $t_{0,n}$ to yield an approximating polynomial $P_{n,t_{0,n}}(t)$ of degree 1 and a remainder of the form $|t-t_{0,n}|^{\alpha}$. Our approximation $x^*_n$ of $T$ is simply the root of the equation $P_{n,t_{0,n}}(t)=1/\eps_n$. Once again, this amounts to dropping the $O_P(|t-t_{0,n}
|^{\alpha})$ from our expansion of $h_n(t)$. The proof ensures that $x^*_n$ has all the properties announced in the Theorem - in particular that it is close to $T$.

\subsubsection{Case $\mysigmasquared_n<\mysigmasquared_{0,n}$}\label{subsubsec:detailedTreatmentPhaseTFirstCase}
We treat this case in full detail and go faster on the two other ones, since the ideas are similar. Recall that the equation defining $T$ is
$$
\frac{1}{\mysigmasquared_n}=h_n(T)=\frac{\sum_{j=1}^{l_n} \coordCanonRandomVector_{j,n}^2}{n}\frac{1}{T}+h_{l_n,n}(T)\;.
%\frac{\sum_{j=1}^{l_n} \coordCanonRandomVector_j^2}{n}\frac{1}{T}+\frac{1}{n}\sum_{j=l+1}^n \frac{\coordCanonRandomVector_j^2}{T+\gamma_n+\delta_j}=
$$

In this case,
$
g_{l_n,n}(0)=\frac{1}{\mysigmasquared_{0,n}}<\frac{1}{\mysigmasquared_n}\;.
$
% so it is clear that the term $\frac{\sum_{j=1}^{l_n} \coordCanonRandomVector_j^2}{n T}$ needs to enter into play to ``saturate" the equality. 
Let us first localize $T$. Denoting $\chi^2_{l_n}=\sum_{j=1}^{l_n} \coordCanonRandomVector_{j,n}^2$, and using  $h_n(t)\geq \chi^2_{l_n}/(nt)$ as well as the fact that $h_n$ is decreasing, we see that $T\geq (\chi^2_{l_n}/n) \mysigmasquared_n$. On the other hand, $h_n(t)\leq \chi^2_{l_n}/(nt)+h_{l_n,n}(0)$. Recall that $h_{l_n,n}(0)=g_{l_n,n}(0)+\gO_P(n^{-1/2})=1/\mysigmasquared_{0,n}+\gO_P(n^{-1/2})$. Simple algebra then gives that $T\leq (\chi^2_{l_n}/n) \mysigmasquared_n/(1-\mysigmasquared_n h_{l_n,n}(0))$. Of course, in the situation we are investigating, $\mysigmasquared_n h_{l_n,n}(0)$ is bounded away from 1 with probability going to 1 as $n\tendsto \infty$.

Let us now expand the last term above, i.e $h_{l_n,n}(t)$, in powers of $t$'s. Because $h'_{l_n,n}$ is uniformly bounded in probability for $t$ in a neighborhood of $0$, we have, for small $t$,
$$
h_{l_n,n}(t)=\frac{1}{n}\sum_{j=l_n+1}^n \frac{\coordCanonRandomVector_{j,n}^2}{\gamma_n+\delta_{j,n}}+\gO_P(t)=\frac{1}{\mysigmasquared_{0,n}}+
\frac{1}{\sqrt{n}}\xi_{1,n}+\gO_P(t)\;,
% \frac{1}{\sqrt{n}}\sum_{j=l_n+1}^n \frac{\coordCanonRandomVector_{j,n}^2-1}{\gamma_n+\delta_{j,n}}+\gO_P(t)
% \;.
% -t \frac{1}{n}\sum_{j=l_n+1}^n \frac{\coordCanonRandomVector_{j,n}^2}{(\gamma_n+\delta_{j,n})^2}+t^2 \frac{1}{n}\sum_{j=l_n+1}^n \frac{\coordCanonRandomVector_{j,n}^2}{(\gamma_n+\delta_{j,n})^3}+\gO_P(t^3)\;.
$$
where $\xi_{1,n}=\frac{1}{\sqrt{n}}\sum_{j=l_n+1}^n \frac{\coordCanonRandomVector_{j,n}^2-1}{\gamma_n+\delta_{j,n}}=\gO_P(1)$.
% So calling $\zeta_1=\frac{1}{n}\sum_{j=l_n+1}^n \frac{\coordCanonRandomVector_{j,n}^2}{(\gamma_n+\delta_{j,n})^2}$, and $\zeta_2=\frac{1}{n}\sum_{j=l_n+1}^n \frac{\coordCanonRandomVector_{j,n}^2}{(\gamma_n+\delta_{j,n})^3}$ the equation defining $T$ becomes
% $$
% \frac{1}{\mysigmasquared}=\frac{\chi^2_{l_n}}{n}\frac{1}{T}+\frac{1}{\mysigmasquared_0}+\frac{1}{n}\sum_{j=l_n+1}^n \frac{\coordCanonRandomVector_{j,n}^2-1}{\gamma_n+\delta_{j,n}}-T \zeta_1+T^2 \zeta_2+\gO(T^3)\;.
% $$
We see that by taking
$$
t(2)=\frac{W_1}{n} \left(1+\frac{1}{\sqrt{n}}\frac{\xi_{1,n}}{\frac{1}{\mysigmasquared_n}-\frac{1}{\mysigmasquared_{0,n}}}\right)\;, 
\text{ with }
W_1=\frac{\chi^2_{l_n}}{\frac{1}{\mysigmasquared_n}-\frac{1}{\mysigmasquared_{0,n}}}\;,
$$
we have
$$
h_n(t(2))-h_n(T)=h_n(t(2))-\frac{1}{\mysigmasquared_n}=\gO_P(1/n)\;.
$$

It is clear that both $t(2)$ and $T$ are contained in the interval $I=\left(\chi^2_{l_n} \mysigmasquared_n /n,2 \mysigmasquared_n \chi^2_{l_n} \upsilon_n/n \right)$, where $\upsilon_n=\max[ 1/(1-\mysigmasquared_n h_{l_n,n}(0)),1/(1-\mysigmasquared_n/\mysigmasquared_{0,n})]$.  The mean value theorem gives 
$$
|t(2)-T|\leq \frac{|h_n(t(2))-h_n(T)|}{\inf_{t\in I} |h_n'(t)|}\;.
$$
Of course, $|h_n'(t)|\geq \chi^2_{l_n}/(nt^2)+|h'_{l_n,n}(0)|\geq \chi^2_{l_n}/(nt^2)$. So we see that $\inf_{t\in I} |h_n'(t)|\geq n R_n $, where $R_n$ is a positive random variable bounded away from 0 with probability going to 1, since we assume that $\mysigmasquared_n$ and $l_n$ remain bounded. We conclude that 
$$
|t(2)-T|=\gO_P(|h_n(t(2))-h_n(T)|/n)=\gO_P(n^{-2})\;, \text{ as announced in Theorem \ref{thm:PhaseTransition}}.
$$
\subsubsection{Case $\mysigmasquared_n=\mysigmasquared_{0,n}$}
We now have, for small $t$, using the fact that $g_{l_n,n}(0)=1/\mysigmasquared_{0,n}=1/\mysigmasquared_n$, 
$$
h_n(t)=\frac{\chi^2_{l_n}}{n}\frac{1}{t}+\frac{1}{\mysigmasquared_n}+\frac{1}{n}\sum_{j=l_n+1}^n \frac{\coordCanonRandomVector_{j,n}^2-1}{\gamma_n+\delta_{j,n}}-t \zeta_{1,n}+\gO_P(t^2)\;,
$$
where $\zeta_{1,n}=\frac{1}{n}\sum_{j=l_n+1}^n \frac{\coordCanonRandomVector_{j,n}^2}{(\gamma_n+\delta_{j,n})^2}=\gO_P(1)$. 
Because $\xi_{1,n}=\frac{1}{\sqrt{n}}\sum_{j=l_n+1}^n \frac{\coordCanonRandomVector_{j,n}^2-1}{\gamma_n+\delta_{j,n}}=\gO_P(1)$, we see that now, $T$ has to be of order $1/\sqrt{n}$, since $h_n(T)=1/\mysigmasquared_n$. Using the ansatz $t(1)=\alpha/\sqrt{n}$, we see that 
$$
h_n(t(1))-\frac{1}{\mysigmasquared_n}=\gO_P\left(\frac{1}{n}\right)\;, \text{ if } \alpha=\frac{\xi_{1,n}+\sqrt{\xi_{1,n}^2+4\chi^2_{l_n} \zeta_{1,n}}}{2 \zeta_{1,n}}\;.
$$
In a neighborhood of $\alpha/\sqrt{n}$, $h_n$ is Lipschitz with Lipschitz constant bounded away from 0, with probability going to 1.
Hence, as argued in \ref{subsubsec:mainArg} and detailed in \ref{subsubsec:detailedTreatmentPhaseTFirstCase}, we can conclude that 
$$
T=\frac{\alpha}{\sqrt{n}}+\gO_P(\frac{1}{n})\;.
$$
\subsubsection{Case $\mysigmasquared_n>\mysigmasquared_{0,n}$}
Recall that the equation defining $T$ is
$$
\frac{1}{\mysigmasquared_n}=\frac{\chi^2_{l_n}}{n}\frac{1}{T}+\frac{1}{n}\sum_{j=l_n+1}^n\frac{\coordCanonRandomVector_{j,n}^2-1}{T+\gamma_n+\delta_{j,n}}+\frac{1}{n}\sum_{j=l_n+1}^n \frac{1}{T+\gamma_n+\delta_{j,n}}\;.
$$
When $\mysigmasquared_n>\mysigmasquared_{0,n}$, we can find $t_{0,n}$ bounded away from 0 such that
$$
\frac{1}{\mysigmasquared_n}=\frac{1}{n}\sum_{j=l_n+1}^n \frac{1}{t_{0,n}+\gamma_n+\delta_{j,n}}\;.
$$
$t_{0,n}$ is furthermore bounded - in $n$ - under our assumptions.

By writing $t=t_{0,n}+\smallQuantity$, for $\smallQuantity$ small, after expanding $h_n$ around $t_{0,n}$, we see that we have
$$
h_n(t)=\frac{\chi^2_{l_n}}{nt_{0,n}}+\frac{1}{\sqrt{n}} \xi(t_{0,n})+\frac{1}{\mysigmasquared_n}-\smallQuantity \zeta(t_{0,n})+\gO_P\left(\max\left[\frac{\smallQuantity}{\sqrt{n}},\smallQuantity^2\right]\right)\;,%=\frac{1}{\mysigmasquared_n}+\gO_P\left(\max\left[\smallQuantity^2,\frac{1}{n},\frac{\eta}{\sqrt{n}} \right]\right)\;,
$$
where
$$
\xi(t_{0,n})=\frac{1}{\sqrt{n}}\sum_{j=l_n+1}^n \frac{\coordCanonRandomVector_{j,n}^2-1}{t_{0,n}+\gamma_n+\delta_{j,n}}=\gO_P(1)\;, \text{and }
\zeta(t_{0,n})=\frac{1}{n}\sum_{j=l_n+1}^n\frac{1}{(t_{0,n}+\gamma_n+\delta_{j,n})^2}=\gO_P(1)\;.
$$

Let us call
$$
t(1)=t_{0,n}+\frac{1}{\sqrt{n}}\frac{\xi(t_{0,n})}{\zeta(t_{0,n})}\;.
$$
Our assumptions and the fact that $t_{0,n}\leq \mysigmasquared_n$ guarantee that $\zeta(t_{0,n})$ is bounded below as $n$ becomes large.
The expansion above shows that
$$
\frac{1}{\mysigmasquared_n}-h_n(t(1))=\gO_P(1/n)\;.
$$
Because, with probability going to 1, $h_n$ is Lipschitz with Lipschitz constant bounded below in a neighborhood of $t_{0,n}$, we conclude as in \ref{subsubsec:mainArg} that
$$
T=t_{0,n}+\frac{1}{\sqrt{n}}\frac{\xi(t_{0,n})}{\zeta(t_{0,n})} +\gO_P\left(\frac{1}{n}\right),
$$
which concludes the proof.
\end{proof}

The phase transition can be further explored in the situation where $\mysigmasquared_n-\mysigmasquared_{0,n}$ is infinitesimal in $n$ but not exactly zero. We are especially concerned in this paper with random variables of the type $$\max_{i=1,\ldots,k} \lambda_{\max}(X_n+(\mysigmasquared_n/n)\canonicalRandomVector_i\canonicalRandomVector_i^{T})-\lambda_{\max}(X_n)$$ for i.i.d $\canonicalRandomVector_i$'s. The previous theorem gives us an idea of the scale of this difference, which clearly depends on $\mysigmasquared_n$ and the whole spectrum of $X_n$. It is also clear that  taking a max over finitely many $k$'s does not change anything to the previous result as far as scale is concerned. The previous theorem shows that our uniform bound on the inverse of the gap cannot be improved: in case (1) of the previous theorem, the gap between the two largest eigenvalues of $X_n(\mysigmasquared_n)$ scales like $1/n$, the rate we obtained in our non-asymptotic bounds. However, in many situations, {\em the gap is much greater than $1/n$}, usually of order at least $1/\sqrt{n}$, and the worst case bound on the Lipschitz constant of $\smoothedF_k(X)$ is very conservative.

\section{Stochastic composite optimization} \label{s:sco}
In this section, we will develop a variant of the algorithm in \citep{Lan09} which allows for adaptive (monotonic) scaling of the step size parameter. For the sake of completeness, we first recall the key definitions in \citep{Lan09}, adopting the same notation, with only a few minor modifications to allow the full problem to be stochastic. We focus on the following optimization problem
\BEQ\label{eq:min-prob}
\min_{x \in Q} \Psi(x) := f(x)+h(x),
\EEQ
where $Q\subset \reals^n$ is a compact convex set. We let $\|\cdot\|$ be a norm and write $\|\cdot\|_*$ the dual norm. We assume that we only observe noisy oracles for $f(x)$ and $h(x)$ written
\[
f(x,\xi) \quad \mbox{and} \quad h(x,\xi),
\]
for some random variable $\xi\in\reals^d$; we write $\Psi(x,\xi):=f(x,\xi)+h(x,\xi)$ with $\Psi(x)=\Expect[\Psi(x,\xi)]$. We also assume that $\Psi(\cdot,\xi)$ is convex for any $\xi\in\reals^d$, and that $\Psi(x,\xi) \geq \Psi(x,0)$ a.s. with 
\[
\Expect[\Psi(x^*,\xi)]-\Psi(x^*,0) \leq \mu
\]
for some $\mu>0$ at the optimum of problem~\eqref{eq:min-prob}, with $\mu$ typically of order $\epsilon$. The value of $\mu$ is typically controlled by the magnitude of the noise $\xi$. The function $f(x)$ is assumed to be convex with Lipschitz continuous gradient
\[
\|\nabla f(x) -\nabla f(y)\|_* \leq L \|x - y\|, \quad \mbox{for all } x,y\in Q,
\]
and $h(x)$ is also assumed to be a convex Lipschitz continuous function with
\[
|h(x)-h(y)| \leq \mathcal{M} \|x-y\|, \quad \mbox{for all } x,y\in Q.
\]
Furthermore, we assume that we observe a subgradient of $\Psi$ through a stochastic oracle $G(x,\xi)$, satisfying
\BEA\label{eq:grad-oracle}
\Expect[G(x,\xi)]=g(x) \in \partial\Psi(x),\\
\Expect[\|G(x,\xi)-g(x)\|_*^2]\leq \sigma^2.
\EEA
We let $\omega(x)$ be a distance generating function, i.e. a function such that
\[
Q^o=\left\{x\in Q:~ \exists y\in \reals^p,~x\in \argmin_{u\in Q} [y^Tu + \omega(u)]\right\}
\]
is a convex set. The function $\omega(x)$ is strongly convex on $Q^o$ with modulus $\alpha$ with respect to the norm $\|\cdot\|$, which means
\[
(y-x)^T(\nabla\omega(y)-\nabla\omega(x)) \geq \alpha \|y-x\|^2, \quad x,y\in Q^o.
\]
We then define a prox-function $V(x,y)$ on $Q^o \times Q$ as follows:
\BEQ\label{eq:prox}
V(x,y)\equiv \omega(y) - [ \omega(x)+\nabla \omega(x)^T(y-x)].
\EEQ
It is nonnegative and strongly convex with modulus $\alpha$ with respect to the norm $\|\cdot\|$. The prox-mapping associated to $V$ is then defined as
\BEQ \label{prox-map}
P_x^{Q,\omega}(y) \equiv \argmin_{z\in Q} \{ y^T(z-x) + V(x,z)\}.
\EEQ
This prox-mapping can be rewritten
\[
P_x^{Q,\omega}(y) = \argmin_{z\in Q} \{ z^T(y-\nabla\omega(x)) + \omega(z) \},
\]
and the strong convexity of $\omega(\cdot)$ means that $P_x^{Q,\omega}(\cdot)$ is Lipschitz continuous with respect to the norm $\|\cdot\|$ with modulus $1/\alpha$ (see \citet{Nemi04} or \cite[Vol.\,II,\,Th.\,4.2.1]{Hiri96}).
Finally, we define the $\omega$ diameter of the set $Q$ as
\BEQ \label{eq:diameter}
D_{\omega,Q}\equiv(\max_{z\in Q} \omega(z)-\min_{z\in Q} \omega(z))^{1/2},
\EEQ
and let
\[
x^\omega=\argmin_{x\in Q} \omega(x),
\]
which satisfies
\[
\frac{\alpha}{2} \|x-x^\omega\|^2 \leq V(x^\omega,x) \leq \omega(x)-\omega(x^\omega) \leq D_{\omega,Q}^2, \quad \mbox{for all }x\in Q.
\]
\citep[Corollary\,1]{Lan09} implies the following result on the complexity of solving~\eqref{eq:smooth-min} using the AC-SA algorithm in \cite[\S3]{Lan09}.

\begin{proposition}\label{prop:cvg-sdp}
Let $N>0$, and write $\smoothedF^*_k$ the optimal value of problem~\eqref{eq:smooth-min}. Suppose that the sequences $X_{t},X_{t}^{md},X_{t}^{ag}$ are computed as in \citep[Corollary\,1]{Lan09} using the stochastic gradient oracle in~\eqref{eq:stoch-oracle}. After $N$ iterations of the AC-SA algorithm in \cite[\S3]{Lan09}, we have
\BEQ\label{eq:expect-prec-bnd}
\Expect[\smoothedF_k(X_{N+1}^{ag})-\smoothedF^*_k] \leq \frac{8n C_k  D_{\omega,Q}^2}{\epsilon N(N+2)} +\frac{4\sqrt{2} D_{\omega,Q}}{\sqrt{Nq}}.
\EEQ
\end{proposition}
\begin{proof}
Using the bound on the variance of the stochastic oracle $G(X,\canonicalRandomVector)$, we know that $G$ satisfies~\eqref{eq:grad-oracle} with $\sigma^2=1/q$. Section~\ref{s:smooth} also shows that the Lipschitz constant of the gradient is bounded by $C_k n/\epsilon$. If we pick $\|\cdot\|_F^2/2$ as the prox function, \citep[Corollary\,1]{Lan09} yields the desired result.
\end{proof}

\noindent Setting $N=2D_{\omega,Q}\sqrt{n}/\epsilon$ and $q=\max\{1,D_{\omega,Q}/(\epsilon \sqrt{n})\}$ in the convergence bound above will then ensure $\Expect[\smoothedF_k(X_N)-\smoothedF_k^*]=O(\epsilon)$. Because our bounds on the Lipschitz constant are usually very conservative, in the section that follows, we detail a version of the AC-SA algorithm with adaptive (but monotonically decreasing) step-size scaling parameter.

\subsection{Stochastic composite optimization with line search} \label{ss:ls}
The algorithm in \citep[\S3]{Lan09} uses worst case values of the Lipschitz constant $L$ and of the gradient's quadratic variation $\sigma^2$ to determine step sizes at each iteration. In practice, this is a conservative strategy and slows down iterations in regions where the function is smoother. In the deterministic case, adaptive versions of the optimal first-order algorithm in \citep{Nest83} have been developed by \citet{Nest07} among others. These algorithms run a few line search steps at each iteration to determine the optimal step size while guaranteeing convergence. The algorithm in \citep{Lan09} is a generalization of the first-order methods in \citep{Nest83,Nest03a} and, in what follows, we adapt the line search steps in \citet{Nest07} to the stochastic algorithm of \citep[\S3]{Lan09}. Here, we will study the convergence properties of an adaptive variant of the algorithm for stochastic composite optimization in \citep[\S3]{Lan09}, with monotonic line search.

\begin{algorithm}[ht]
\caption{Adaptive algorithm for stochastic composite optimization.}
\label{alg:stoch}
\begin{algorithmic} [1]
\REQUIRE An initial point $x^{ag}=x_1=x^w \in\reals^n$, an iteration counter $t=1$, the number of iterations $N$, line search parameters $\gamma^{min},\gamma^{max},\,\gamma^{d},\,\gamma>0$, with $\gamma^d<1$.
\STATE Set $\gamma=\gamma^{max}$.
\FOR{$t=1$ to $N$}
\STATE Define $x_t^{md}=\frac{2}{t+1} x_t + \frac{t-1}{t+1} x_t^{ag}$
\STATE Call the stochastic gradient oracle to get $G(x_t^{md},\xi_t)$.
\REPEAT
\STATE Set $\gamma_t=\frac{(t+1)\gamma}{2}$.
\STATE Compute the prox mapping $x_{t+1}=P_{x_t}(\gamma_t G(x_t^{md},\xi_t))$.
\STATE Set $x_{t+1}^{ag}=\frac{2}{t+1} x_{t+1} + \frac{t-1}{t+1} x_t^{ag}$.
\UNTIL{$\Psi(x_{t+1}^{ag},\xi_{t+1})\leq\Psi(x_{t}^{md},\xi_t) + \langle G(x_t^{md},\xi_t), x_{t+1}^{ag}-x_{t}^{md}\rangle+ \frac{\alpha\gamma^d}{4 \gamma} \|x_{t+1}^{ag}-x_{t}^{md}\|^2 + 2 \mathcal{M} \|x_{t+1}^{ag}-x_{t}^{md}\|$ or $\gamma\leq \gamma^{min}$. If exit condition fails, set $\gamma=\gamma\gamma^d$ and go back to step 5.}
\STATE Set $\gamma=\max\left\{\gamma^{min},\gamma \right\}$.
\ENDFOR
\ENSURE A point $x^{ag}_{N+1}$.
\end{algorithmic}
\end{algorithm}

In this section, we first modify the convergence lemma in \citep[Lemma 5]{Lan09} to adapt it to the line search strategy detailed in Algorithm~\ref{alg:stoch}. Note that our method requires testing the line search exit condition using {\em two} oracle calls, the current one in $\xi_t$ and the next one in $\xi_{t+1}$. This last oracle call is of course recycled at the next iteration.

\begin{lemma}\label{lem:step}
Assume that $\Psi(\cdot,\xi_t)$ is convex for any given sample  of the r.v. $\xi_t$. Let $x_{t},x_{t}^{md},x_{t}^{ag}$ be computed as in Algorithm~\ref{alg:stoch}, with $\beta_t=(t+1)/2$. Suppose also that $\gamma$ and these points satisfy the line search exit condition in line 9, i.e.
\[
\Psi(x_{t+1}^{ag},\xi_{t+1})\leq\Psi(x_{t}^{md},\xi_t) + \langle G(x_t^{md},\xi_t), x_{t+1}^{ag}-x_{t}^{md}\rangle+ \frac{\alpha \beta_t}{4 \gamma_t } \|x_{t+1}^{ag}-x_{t}^{md}\|^2 + 2 \mathcal{M} \|x_{t+1}^{ag}-x_{t}^{md}\|
\]
then, for every $x$ in the feasible set, we have
\BEAS
\beta_t\gamma_t[\Psi(x_{t+1}^{ag},\xi_{t+1})-\Psi(x,0)] + V(x_{t+1},x) & \leq & (\beta_t-1)\gamma_t [\Psi(x_{t}^{ag},\xi_t)-\Psi(x,0)] + V(x_t,x)\\
&& + \gamma_t (\Psi(x,\xi_t)-\Psi(x,0))+ \frac{4\mathcal{M}^2\gamma_t^2}{\alpha}.
\EEAS
\end{lemma}
\begin{proof}
As in \citep[Lemma 5]{Lan09}, we write $d_t=x_{t+1}-x_{t}$ and use the parameter $\beta_t=(t+1)/2$ for step sizes so that $x_{t+1}^{ag}-x_{t}^{md}=d_t/\beta_t$. If the current iterates satisfy the line search exit condition, the fact that $\alpha\|d_t\|^2/2\leq V(x_t,x_{t+1})$ by construction yields
\BEAS
\beta_t\gamma_t\Psi(x_{t+1}^{ag},\xi_{t+1}) & \leq & \beta_t\gamma_t[\Psi(x_{t}^{md},\xi_t) + \langle G(x_t^{md},\xi_t), x_{t+1}^{ag}-x_{t}^{md}\rangle] + \frac{\alpha}{4} \|d_t\|^2 + 2\gamma_t \mathcal{M} \|d_t\|\\
& \leq & \beta_t\gamma_t[\Psi(x_{t}^{md},\xi_t) + \langle G(x_t^{md},\xi_t), x_{t+1}^{ag}-x_{t}^{md}\rangle] + V(x_t,x_{t+1}) - \frac{\alpha}{4} \|d_t\|^2 + 2\gamma_t \mathcal{M} \|d_t\|.
\EEAS
Using the convexity of $\Psi(\cdot,\xi_t)$ we then get 
\BEAS
&& \beta_t\gamma_t[\Psi(x_{t}^{md},\xi_t) + \langle G(x_t^{md},\xi_t), x_{t+1}^{ag}-x_{t}^{md}\rangle] \\
& = & (\beta_t-1)\gamma_t[\Psi(x_{t}^{md},\xi_t) + \langle G(x_t^{md},\xi_t), x_{t}^{ag}-x_{t}^{md}\rangle] + \gamma_t[\Psi(x_{t}^{md},\xi_t) + \langle G(x_t^{md},\xi_t), x_{t+1}-x_{t}^{md}\rangle]\\
& \leq & (\beta_t-1)\gamma_t \Psi(x_{t}^{ag},\xi_t) + \gamma_t[\Psi(x_{t}^{md},\xi_t) + \langle G(x_t^{md},\xi_t), x_{t+1}-x_{t}^{md}\rangle].
\EEAS
Combining these last two results and using the fact that $bu - au^2/2 \leq b^2/(2a)$ whenever $a>0$, we obtain
\BEAS
\beta_t\gamma_t\Psi(x_{t+1}^{ag},\xi_{t+1}) & \leq & (\beta_t-1)\gamma_t \Psi(x_{t}^{ag},\xi_t) +\gamma_t[\Psi(x_{t}^{md},\xi_t) + \langle G(x_t^{md},\xi_t), x_{t+1}-x_{t}^{md}\rangle] \\
& & + V(x_t,x_{t+1}) - \frac{\alpha}{4} \|d_t\|^2 + 2\gamma_t \mathcal{M} \|d_t\|\\
& \leq & (\beta_t-1)\gamma_t \Psi(x_{t}^{ag},\xi_t) +\gamma_t[\Psi(x_{t}^{md},\xi_t) + \langle G(x_t^{md},\xi_t), x_{t+1}-x_{t}^{md}\rangle] \\
& & + V(x_t,x_{t+1}) + \frac{4\gamma_t^2\mathcal{M}^2}{\alpha}.
\EEAS
For any $x$ in the feasible set, we can then use the properties of the prox mapping detailed in \citep[Lemma 1]{Lan09}, with $p(\cdot)=\gamma_t \langle G(x_t^{md},\xi_t), \cdot-x_{t}^{md}\rangle$ together with the convexity of $\Psi(\cdot,\xi_t)$ and the definition of $x_{t+1}$ in Algorithm~\ref{alg:stoch} to show that
\BEAS
&& \gamma_t[\Psi(x_{t}^{md},\xi_t) + \langle G(x_t^{md},\xi_t), x_{t+1}-x_{t}^{md}\rangle] + V(x_t,x_{t+1}) \\
& \leq & \gamma_t \Psi(x_{t}^{md},\xi_t) + \gamma_t \langle G(x_t^{md},\xi_t), x-x_{t}^{md}\rangle + V(x_t,x) - V(x_{t+1},x)\\
& \leq & \gamma_t \Psi(x,\xi_t) + V(x_t,x) - V(x_{t+1},x)\;.
\EEAS
Combining these last results shows that
\[
\beta_t\gamma_t\Psi(x_{t+1}^{ag},\xi_{t+1}) \leq (\beta_t-1)\gamma_t \Psi(x_{t}^{ag},\xi_t) + \gamma_t \Psi(x,\xi_t) + V(x_t,x) - V(x_{t+1},x) + \frac{4\gamma_t^2\mathcal{M}^2}{\alpha}\;,
\]
and subtracting $\beta_t\gamma_t\Psi(x,0)$ from both sides yields the desired result.
\end{proof}

We are now ready to prove the main convergence result, adapted from \cite[Corollary 1]{Lan09}. We simply stitch together the convergence results we obtained in Lemma~\ref{lem:step} for the line search phase of the algorithm, with that of \citep[Lemma 5]{Lan09} for the second phase where $\gamma=\gamma^{min}$, writing the switch time~$T_\gamma$. Note that the step size is still increasing in the second phase of the algorithm because $\gamma_t=\gamma^{min}(t+1)/2$.

\begin{proposition}\label{prop:expect-convergence}
Let $N>0$, and write $\Psi(x^{*},0)$ the optimal value of problem~\eqref{eq:min-prob}. Suppose that the sequences $x_{t},x_{t}^{md},x_{t}^{ag}$ are computed as in Algorithm~\ref{alg:stoch}, with line search parameter $\gamma$ initially set to $\gamma=\gamma^{max}$ with
\BEQ\label{eq:gamma-max}
\gamma^{max}\leq \frac{\sqrt{6\alpha} D_{\omega,Q}}{(N+2)^{3/2}(4\mathcal{M}^2+\sigma^2)^{1/2}} \quad \mbox{and} \quad  \gamma^{min}=\min\left\{\frac{\alpha}{2L},\gamma^{max}\right\}
\EEQ
with $\gamma^{d}<1$. After $N$ iterations of Algorithm~\ref{alg:stoch}, we have
\BEQ%\label{eq:expect-prec-bnd}
\Expect[\Psi(x_{N+1}^{ag})-\Psi(x^{*},0)] \leq \frac{8LD_{\omega,Q}^2}{\alpha N^2} + \frac{8}{N^2\gamma^{min}}\Expect\left[ \frac{2(4\mathcal{M}^2+\sigma^2)}{\alpha}\sum_{t=1}^{N}\gamma_t^2 \right]
+ \frac{(T_\gamma+2)^2\gamma^{max}\mu}{N^2 2\gamma^{min}}
\EEQ
and a simpler, but coarser bound is given by
\BEQ\label{eq:expect-coarse-bnd}
\Expect\left[\Psi(x_{N+1}^{ag})-\Psi(x^{*},0)\right]\leq \frac{8LD_{\omega,Q}^2}{\alpha N^2} + \frac{8D_{\omega,Q}\sqrt{4 \mathcal{M}^2 + \sigma^2}}{\sqrt{N}} \left(\frac{\gamma^{max}}{\gamma^{min}}\rho_N+1-\rho_N\right) + \frac{(T_\gamma+2)^2\gamma^{max}\mu}{N^2 2\gamma^{min}},
\EEQ
where $\rho_N=(T_\gamma+2)^3/(N+2)^3$.
\end{proposition}
\begin{proof}
Lemma~\ref{lem:step} applied at $x^*$ shows
\BEAS
\beta_t\gamma_t[\Psi(x_{t+1}^{ag},\xi_{t+1})-\Psi(x^*,0)] + V(x_{t+1},x^*) & \leq & (\beta_t-1)\gamma_t [\Psi(x_{t}^{ag},\xi_t)-\Psi(x^*,0)] + V(x_t,x^*)\\
&& + \frac{4\mathcal{M}^2\gamma_t^2}{\alpha} + \gamma_t (\Psi(x^*,\xi_t)-\Psi(x^*,0))
\EEAS
hence, having assumed $\Psi(x,\xi_t)-\Psi(x,0) \geq 0$ a.s.,

\BEAS
(\beta_{t+1}-1)\gamma_{t}[\Psi(x_{t+1}^{ag},\xi_{t+1})-\Psi(x^*,0)] &\leq& \beta_{t}\gamma_t[\Psi(x_{t+1}^{ag},\xi_{t+1})-\Psi(x^*,0)]\\
&\leq& (\beta_t-1)\gamma_t [\Psi(x_{t}^{ag},\xi_{t})-\Psi(x^*,0)] + \frac{4\mathcal{M}^2\gamma_t^2}{\alpha}\\
&& + \gamma_t (\Psi(x^*,\xi_t)-\Psi(x^*,0)) + V(x_t,x^*) - V(x_{t+1},x^*)
\EEAS
whenever the line search successfully terminates, with the last term satisfying
\[
\Expect[\gamma_t (\Psi(x^*,\xi_t)-\Psi(x^*,0))] \leq \frac{\gamma^{max}(t+1)}{2} \Expect[\Psi(x^*,\xi_t)-\Psi(x^*,0)]\leq \frac{\gamma^{max}(t+1)}{2} \mu
\]
using again $\Psi(x^*,\xi_t)-\Psi(x^*,0) \geq 0$ a.s. When the line search fails, $\gamma_t=\gamma^{min}(t+1)/2$ is deterministic and \citep[Lem.\,5 \& Th.\,2]{Lan09} show that
\[
(\beta_{t+1}-1)\gamma_t[\Psi(x_{t+1}^{ag})-\Psi(x^*,0)]
\leq (\beta_t-1)\gamma_t [\Psi(x_{t}^{ag})-\Psi(x^*,0)] + V(x_t,x^*)- V(x_{t+1},x^*) + \Delta(x^*)
\]
where
\[
\Delta(x^*)\leq \gamma_t\langle \delta_t , x^*-x_t \rangle +\frac{2(4\mathcal{M}^2+\|\delta_t \|_*^2)\gamma_t^2}{\alpha}.
\]
with $\delta_t=G(x_{t}^{md},\xi_t)-g(x_{t}^{md})$ and $\gamma_t\langle \delta_t , x^*-x_t \rangle \leq \gamma_t \|\delta_t\|_*\|x^*-x_t\|$. We call $t=T_{\gamma}+1$ the iteration where the line search first fails. Combining these last results, using $\beta_1=1$, we obtain
\BEAS
&& (\beta_{N+1}-1)\gamma_N \Expect[\Psi(x_{N+1}^{ag})-\Psi(x^*,0)]\\
 & \leq & D_{\omega,Q}^2 + \sum_{t=1}^{T_\gamma} \Expect\left[\frac{4\mathcal{M}^2\gamma_t^2}{\alpha} \right]+ \sum_{t=1}^{T_\gamma} \frac{\gamma^{max}(t+1)}{2} \mu + (\beta_{T_\gamma+1}-1)\gamma_{T_\gamma} \Expect[\Psi(x_{T_\gamma+1}^{ag})-\Psi(x_{T_\gamma+1}^{ag},\xi_{T_\gamma+1})]\\
&& + \sum_{T_\gamma+1}^{N} \Expect\left[\gamma_t\langle \delta_t , x^*-x_t \rangle +\frac{2(4\mathcal{M}^2+\|\delta_t \|_*^2)\gamma_t^2}{\alpha}\right]\\
& \leq & D_{\omega,Q}^2 + \sum_{t=1}^{T_\gamma} \Expect\left[\frac{4\mathcal{M}^2\gamma_t^2}{\alpha}\right] + \displaystyle \sum_{T_\gamma+1}^{N} \Expect\left[\frac{2(4\mathcal{M}^2+\|\delta_t \|_*^2)\gamma_t^2}{\alpha}\right]+ \frac{(T_\gamma+2)^2\gamma^{max}\mu}{4}\\
& \leq & D_{\omega,Q}^2 +  \displaystyle \Expect\left[ \frac{2(4\mathcal{M}^2+\sigma^2)}{\alpha}\sum_{t=1}^{N}\gamma_t^2\right]+ \frac{(T_\gamma+2)^2\gamma^{max}\mu}{4}
\EEAS
because $\Expect[\Psi(x_{T_\gamma}^{ag})-\Psi(x_{T_\gamma}^{ag},\xi_t)]=0$. Using the fact that $\sum_{t=1}^{N} (t+1)^q \leq (N+2)^{q+1}/(q+1)$ for $q=1,2$ then yields the coarser bound.
\end{proof}

We observe that, as in \citep{Nest07}, allowing a line search slightly increases the complexity bound, by a factor
\[
\left(\frac{\gamma^{max}}{\gamma^{min}}\rho(T_\gamma,N)+1-\rho(T_\gamma,N)\right),
\]
where $\rho(T_\gamma,N)=(T_\gamma+2)^3/(N+2)^3$. We will see however that overall numerical performance can significantly improve because the algorithm takes longer steps.

\subsection{Stochastic composite optimization for semidefinite optimization} \label{ss:stoch-sdp}
We can use the results above to solve problem~\eqref{eq:smooth-min}. In this case,
\[
\Psi(X)=\Expect\left[\Psi(X,z)\right]=\Expect\left[\max_{i=1,\ldots,k} \lambdamax\left(X+\frac{\epsilon}{n}\canonicalRandomVector_i\canonicalRandomVector_{i}^T\right)\right]
\]
and by construction $\Psi(X,z) \geq \Psi(X,0)=\lambdamax(X)$. Recall, that with this choice of oracle, Section 2 shows
\[
\lambdamax(X) \leq \Expect\left[\max_{i=1,\ldots,k} \lambdamax\left(X+\frac{\epsilon}{n}\canonicalRandomVector_i\canonicalRandomVector_{i}^T\right)\right] \leq \lambdamax(X) + k \epsilon
\]
so $\mu=k\epsilon$ in Proposition~\ref{prop:expect-convergence}. We use the following gradient oracle
\BEQ\label{eq:stoch-oracle}
G(X,\canonicalRandomVector)=\frac{1}{q}\sum_{l=1}^q \phi_l\phi_l^T
\EEQ
where each $\phi_l$ is a leading eigenvector of the matrix $X+\frac{\epsilon}{n}\canonicalRandomVector_{i_0}\canonicalRandomVector_{i_0}^{T}$, with
$$
i_0=\argmax_{i=1,\ldots,k} \lambdamax\left(X+\frac{\epsilon}{n}\canonicalRandomVector_i\canonicalRandomVector_i^{T}\right),
$$
where $\canonicalRandomVector_i$ are i.i.d. Gaussian vectors $\canonicalRandomVector_i\sim {\mathcal N}(0, {\idm_n} )$ and $k>0$ is a small constant (typically 3) and $q$ is used to control the variance. The Lipschitz constant of the gradient is bounded by~\eqref{eq:Lip} with
\[
L \leq  \frac{n}{\epsilon}\frac{k}{(k-2)}
\]
and the variance of the gradient oracle is bounded by $1/q$ with $\mathcal{M}=0$ in the results above.

\section{Numerical Experiments} \label{s:numexp}
We test the algorithm detailed above on a maximum eigenvalue minimization problem over a hypercube, a problem used in approximating sparse eigenvectors \citep{dAsp04a}. We seek to solve
\BEQ\label{eq:dspca}
\BA{ll}
\mbox{minimize} & \lambdamax(A+X)\\
\mbox{subject to} & -\rho \leq X_{ij} \leq \rho, \quad \mbox{for }i,j=1,\ldots,n
\EA\EEQ
which is a semidefinite program in the matrix $X\in\symm_n$. Since randomly generated matrices $A$ have a highly structured spectrum, we use a covariance matrix from the gene expression data set in \citep{Alon99} to generate the matrix $A\in\symm_n$, varying the number of genes to change the problem dimension $n$ (we select the $n$ genes with the highest variance) and normalizing the matrix $A$ so that its spectral norm is one. We set $\rho=\max\{\diag(A)\}/2$ in \eqref{eq:dspca}.

We also test performance on the classical {\em MaxCut} relaxation. The primal semidefinite program is written
\[
\BA{ll}
\mbox{maximize} & \Tr(CX)\\
\mbox{subject to} & \diag(X)=1,\,X\succeq 0,
\EA\]
in the variable $X\in\symm_n$. The objective matrix $C$ is sampled from the Wishart distribution with $C=G^TG/\|G\|_2^2$ where $G$ is a standard Gaussian matrix. Here, we solve the dual, written
\BEQ\label{eq:maxcut}
\min ~ \lambdamax(C+\diag(w))-\ones^Tw
\EEQ
in the variable $w\in\reals^n$. The problem is unconstrained, and we add a bound on the Euclidean norm of the vector $w$. The prox function used in both examples (where the feasible sets are an hypercube and an Euclidean ball) is the square Euclidean norm, which means that the prox is simply an Euclidean projection of the matrix $X$ in~\eqref{eq:dspca} (with the projection taken elementwise) and of the vector $w$ in~\eqref{eq:maxcut}.

We first compare the performance of Algorithm~\ref{alg:stoch} with that of the corresponding deterministic algorithms, ACSA as detailed in \cite{Lan09} and the accelerated first-order method (with line search) in \cite[\S4]{Nest07} after smoothing problem~\eqref{eq:dspca} as in \citep{Nest04a,dAsp04a}. We set a fixed number of outer iterations for Algorithm~\ref{alg:stoch} and record the number of iterations (and eigenvector evaluations, these numbers differ because of line search steps) required by the algorithm in \cite[\S4]{Nest07} to reach the best objective value attained by the stochastic method. We set $\epsilon=5\times 10^{-2}$, $q=0.1/\epsilon$, $k=3$ and the maximum number of iterations to $O(\sqrt{n})$ in the stochastic algorithm. In line with the discussion of Section~\ref{ss:phase-trans}, we scale down the Lipschitz constant by a factor 100 in both stochastic and deterministic algorithms. This significantly speeds up the algorithms with no apparent effect on convergence, thus confirming that the worst case bounds are indeed somewhat conservative. 

\begin{figure}[ht!]
\begin{center}
\begin{tabular}{cc}
\psfrag{eigvs}[t][b]{Number of Eigs.}
\psfrag{obj}[b][t]{Objective}
\includegraphics[width=0.45 \textwidth]{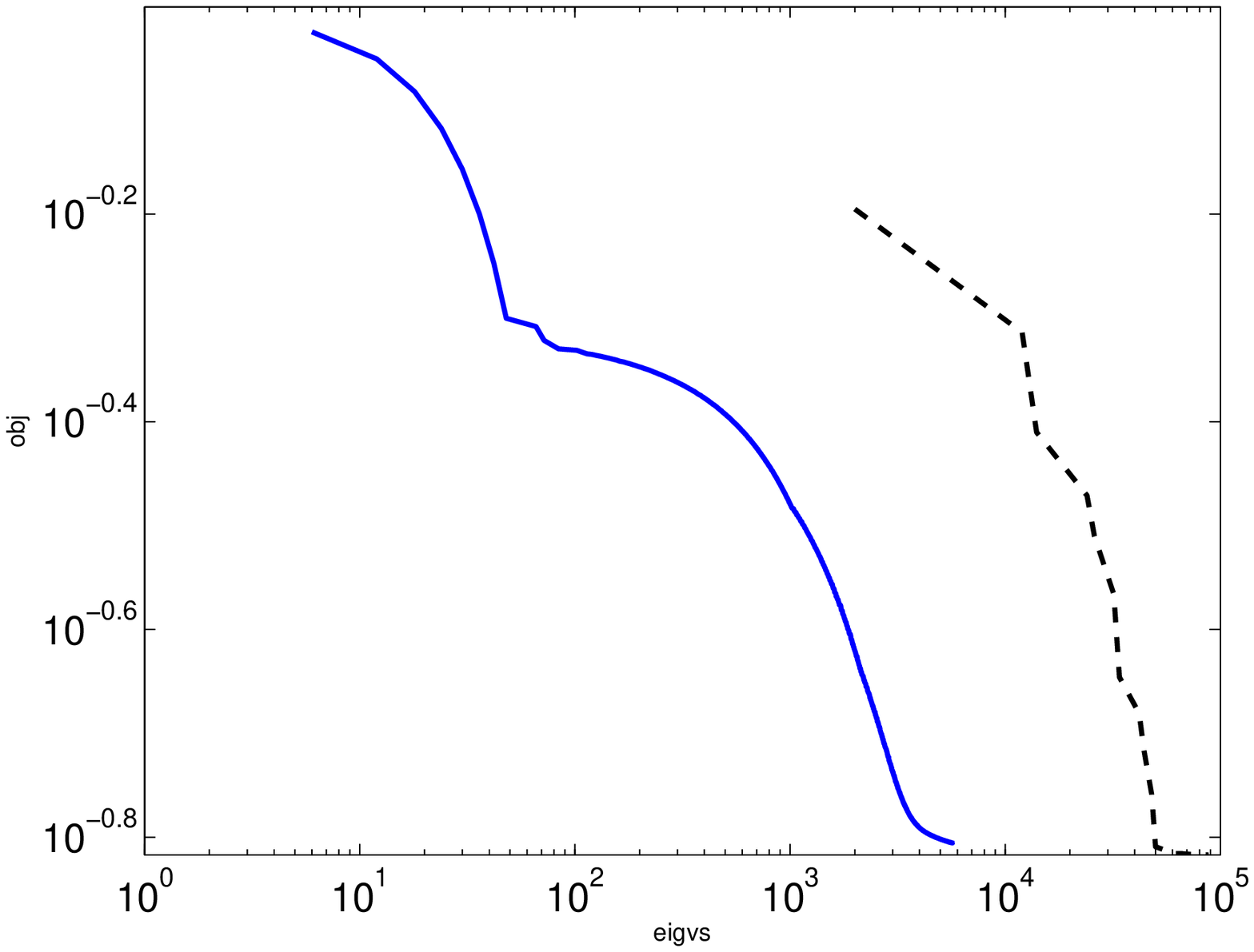}&
\psfrag{eigvs}[t][b]{Number of Eigs.}
\psfrag{obj}[b][t]{Gap}
\includegraphics[width=0.45 \textwidth]{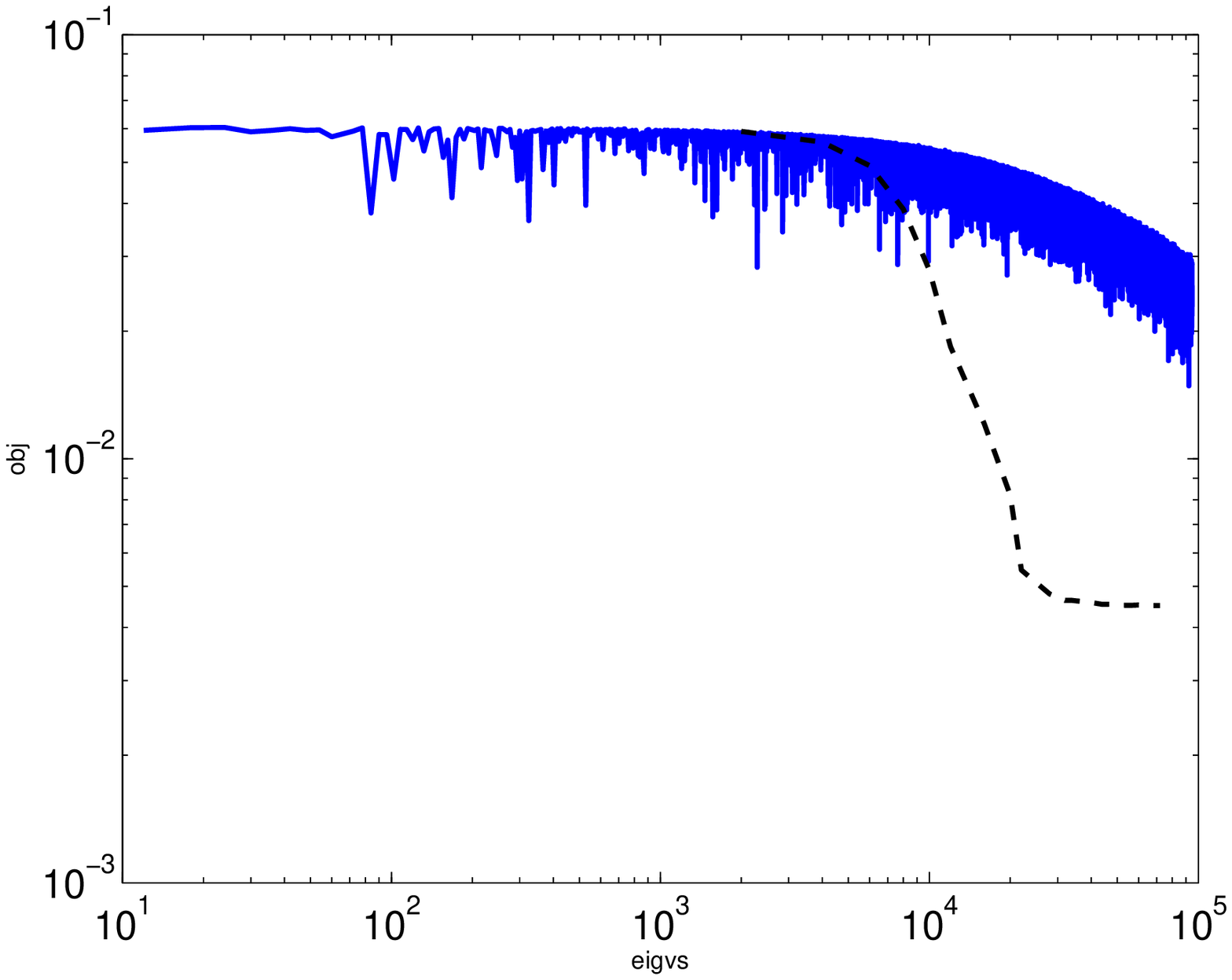}
\end{tabular}
\caption{Objective value (or gap) versus number of eigenvectors computed by Algorithm~\ref{alg:stoch} (solid blue) and the algorithm in \cite[\S4]{Nest07} (dashed black) for the DSPCA {\em (left)} and {\em MaxCut} {\em (right)} relaxations.\label{fig:path-eigs}}
\end{center}
\end{figure}

\begin{figure}[ht!]
\begin{center}
\begin{tabular}{cc}
\psfrag{lam}[t][b]{\lambda_i}
\includegraphics[width=0.45 \textwidth]{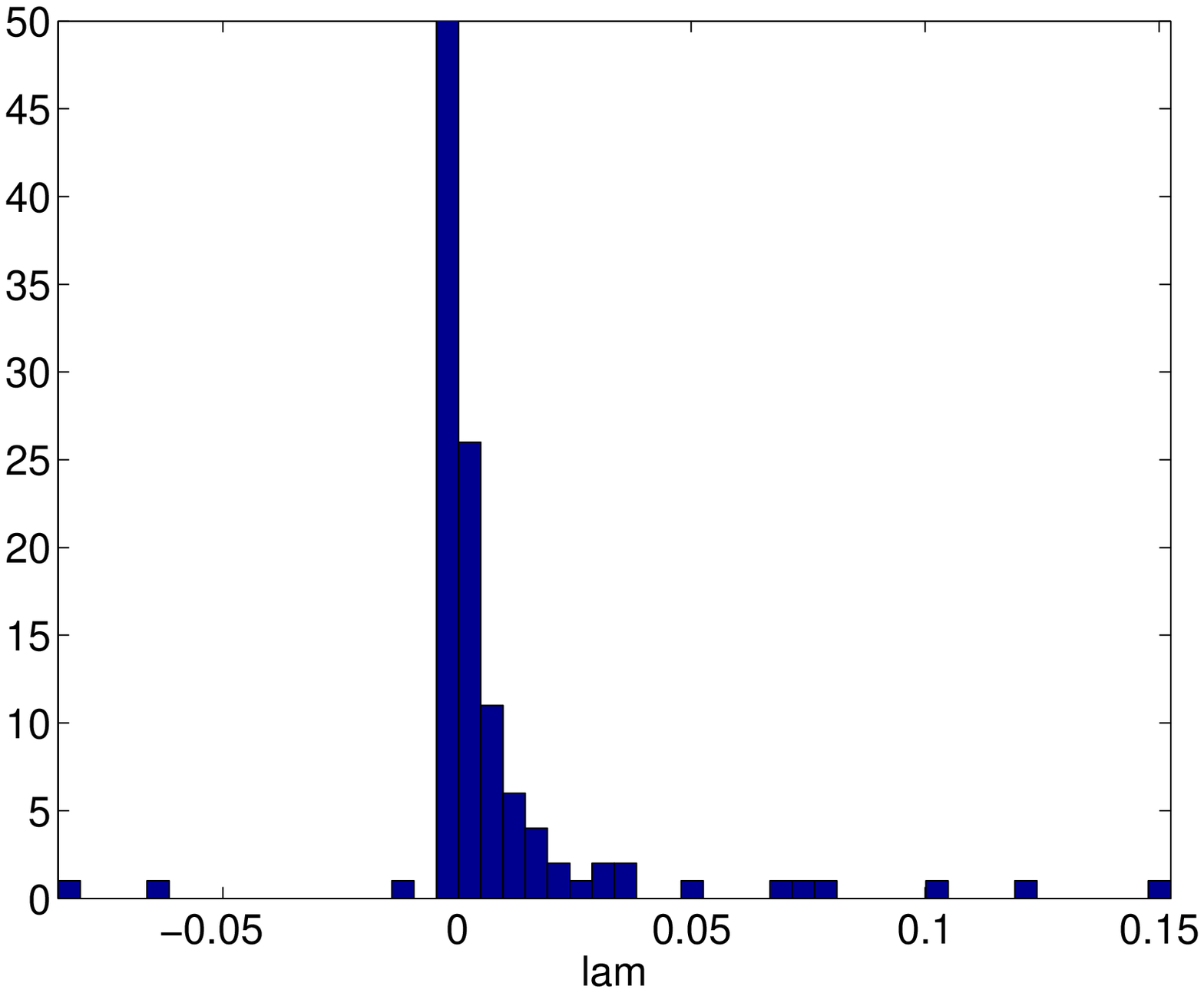}&
\psfrag{lam}[t][b]{\lambda_i}
\includegraphics[width=0.45 \textwidth]{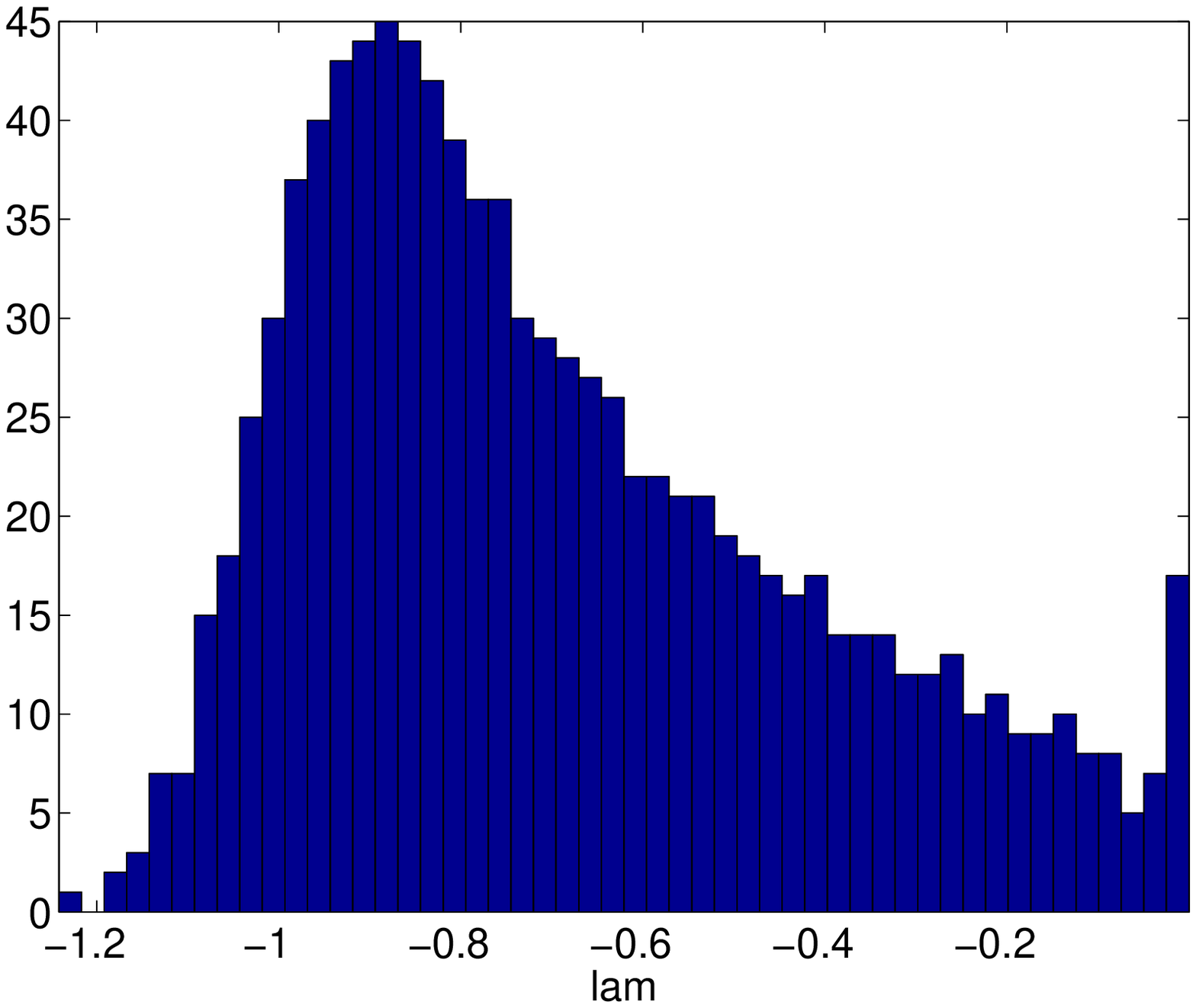}
\end{tabular}
\caption{Histogram of eigenvalues for the matrix solutions to the sparse PCA {\em (left)} and {\em MaxCut} {\em (right)} problems for $n=1000$. For clarity, the graph on the left is truncated above 50.\label{fig:spectrum}}
\end{center}
\end{figure}

To provide a complexity benchmark that is both hardware and implementation independent, we record the total number of eigenvectors used by each algorithm to reach a given objective value (the matrix exponential thus counts as $n$ eigenvectors). We report these results in Tables~\ref{tab:iters} and~\ref{tab:iters-mc} for DSPCA~\eqref{eq:dspca} and {\em MaxCut} respectively. We observe that, for the DSPCA tests, the total number of eigenvectors computed is significantly lower, while the number of iterations is much higher for the stochastic code. The tradeoff is much less favorable for the {\em MaxCut} experiments. In Figure~\ref{fig:path-eigs} we plot the objective value reached as a function of the number of eigenvectors computed for both experiments, when $n=1000$. We again see that the behavior of the stochastic algorithm is much better for DSPCA than for {\em MaxCut}. In Figure~\ref{fig:spectrum}, we plot the spectrum of the solution matrices for both problems. We notice that the leading eigenvalues are much more separated in the DSPCA problem which at least partly explains the difference in performance. More importantly, the deterministic implementation of the ACSA algorithm in \citep{Lan09} seems to be significantly slower than that of the smooth algorithm in \citep{Nest04a}. Improving the numerical performance of the ACSA algorithm itself thus seems to be the key to a competitive implementation of the results detailed here.

%\begin{table}[H]
%\begin{center}
%%\extrarowheight 2ex
%\begin{tabular}{r|c|c|c|c}
%$n$ & \# Iters. (Stoch.) & \# Eigvs. (Stoch.) & \# Iters. (Det.) & \# Eigvs. (Det.) \\
%\hline
%50 & 707 & 1686 & 20 & 4400 \\
%100 & 1000 & 1806 & 20 & 8800 \\
%200 & 1414 & 5052 & 40 & 33600 \\
%500 & 2236 & 6684 & 20 & 45000 \\
%1000 & 3162 & 11406 & 20 & 88000  \end{tabular}
%\vskip 0.5ex
%\caption{Number of iterations and total number of eigenvectors computed by Algorithm~\ref{alg:stoch} (Stoch.) and the algorithm in \cite[\S4]{Nest07} (Det.) to reach identical objective values when solving the DSPCA relaxation in~\eqref{eq:dspca}.\label{tab:iters}}
%\end{center}
%\end{table}

\begin{table}[ht]
\begin{center}
%\extrarowheight 2ex
\begin{tabular}{r|c|c|c|c|c|c}
 & Stoch. & Stoch. & ACSA & ACSA & Det. & Det. \\
$n$ & \# iters. & \# eigvs. & \# iters. & \# eigvs. & \# iters. & \# eigvs.\\
\hline
50 & 707 & 1266 & 51 & 2550 & 16 & 3700 \\
100 & 1000 & 1806 & 50 & 5000 & 12 & 5800 \\
200 & 1414 & 2532 & 55 & 11000 & 28 & 24800 \\
500 & 2236 & 8016 & 60 & 30000 & 12 & 29000 \\
1000 & 3162 & 18990 & 65 & 65000 & 12 & 56000 \\
2000 & 4472 & 21444 & 66 & 132000 & 14 & 132000   
\end{tabular}
\vskip 0.5ex
\caption{Number of iterations and total number of eigenvectors computed by Algorithm~\ref{alg:stoch} (Stoch.), the ACSA algorithm in \cite{Lan09} and the algorithm in \cite[\S4]{Nest07} (Det.) (both with exponential smoohting) to reach identical objective values when solving the DSPCA relaxation in~\eqref{eq:dspca}.\label{tab:iters}}
\end{center}
\end{table}

%\begin{table}[H]
%\begin{center}
%%\extrarowheight 2ex
%\begin{tabular}{r|c|c|c|c}
%$n$ & \# Iters. (Stoch.) & \# Eigvs. (Stoch.) & \# Iters. (Det.) & \# Eigvs. (Det.) \\
%\hline
%50 & 3536 & 21240 & 20 & 4100\\
%100 & 5000 & 30024 & 20 & 8000\\
%200 & 7071 & 42438 & 20 & 15600\\
%500 & 11180 & 67086 & 20 & 37000\\
%1000 & 15811 & 94872 & 20 & 72000   
%\end{tabular}
%\vskip 0.5ex
%\caption{Number of iterations and total number of eigenvectors computed by Algorithm~\ref{alg:stoch} (Stoch.) and the algorithm in \cite[\S4]{Nest07} (Det.) to reach identical objective values when solving the {\em MaxCut} relaxation in~\eqref{eq:maxcut}.\label{tab:iters-mc}}
%\end{center}
%\end{table}

\begin{table}[ht]
\begin{center}
%\extrarowheight 2ex
\begin{tabular}{r|c|c|c|c|c|c}
 & Stoch. & Stoch. & ACSA & ACSA & Det. & Det. \\
$n$ & \# iters. & \# eigvs. & \# iters. & \# eigvs. & \# iters. & \# eigvs.\\
\hline
50 & 3536 & 9534 & 217 & 10850 & 2 & 400 \\
100 & 5000 & 30024 & 353 & 35300 & 4 & 1600 \\
200 & 7071 & 42438 & 537 & 107400 & 6 & 4400 \\
500 & 11180 & 67086 & 545 & 272500 & 6 & 9000 \\
1000 & 15811 & 94872 & 601 & 601000 & 6 & 16000 \\
2000 & 22361 & 134178 & 377 & 754000 & 4 & 20000   
\end{tabular}
\vskip 0.5ex
\caption{Number of iterations and total number of eigenvectors computed by Algorithm~\ref{alg:stoch} (Stoch.), the ACSA algorithm in \cite{Lan09} and the algorithm in \cite[\S4]{Nest07} (Det.) (both with exponential smoohting) to reach identical objective values when solving the {\em MaxCut} relaxation in~\eqref{eq:maxcut}.\label{tab:iters-mc}}
\end{center}
\end{table}

In both algorithms, the cost of each iteration is dominated by that of computing gradients. The cost of each gradient computation in Algorithm~\ref{alg:stoch} is dominated by the cost of computing the leading eigenvector of $q$ perturbed matrices, which is $O(qn^2\log n)$. The cost of each gradient computation in \cite[\S4]{Nest07} is dominated by the cost of computing a matrix exponential, which is $O(n^3)$. This means that the ratio between these costs grows as $O(n/(q\log n))$.

\begin{figure}[ht!]
\begin{center}
\begin{tabular}{cc}
\psfrag{iter}[t][b]{Iteration}
\psfrag{gamma}[b][t]{$\gamma$}
\includegraphics[width=0.45 \textwidth]{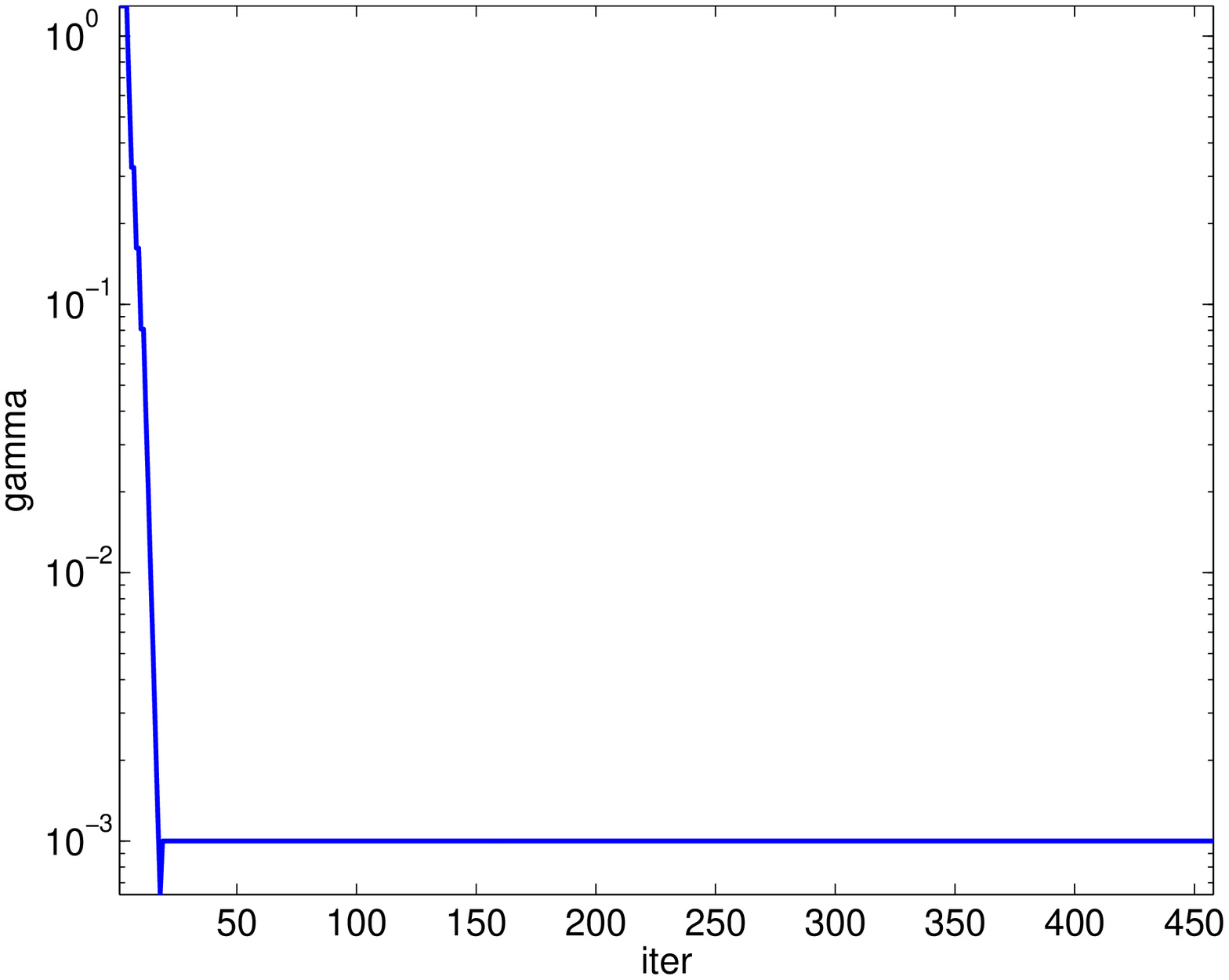}&
\psfrag{iter}[t][b]{Iteration}
\psfrag{L}[b][t]{$1/L$}
\includegraphics[width=0.45 \textwidth]{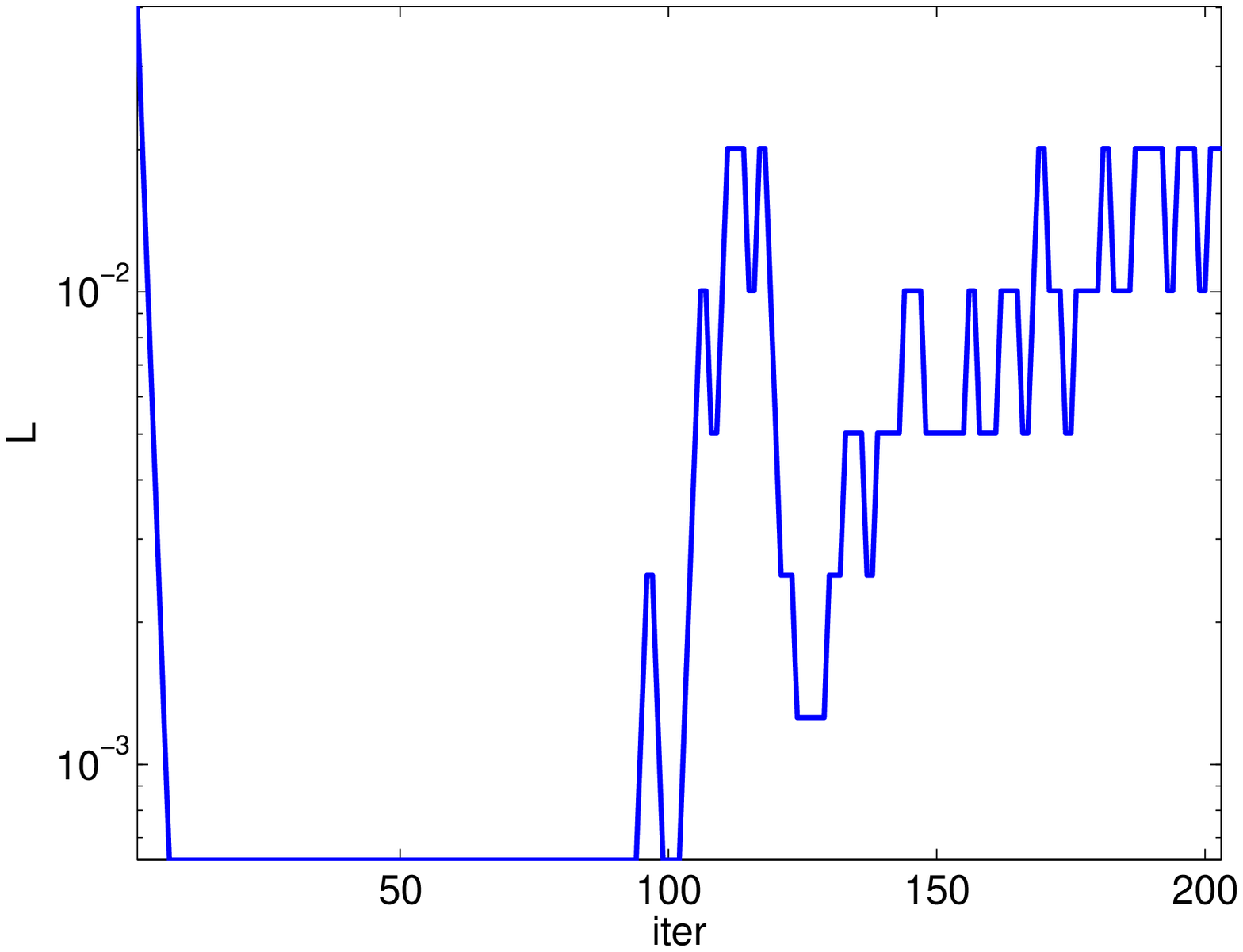}
\end{tabular}
\caption{Line search parameters $\gamma$ for the stochastic algorithm {\em (left)} together with the values of the inverse of the Lipschitz constant $L$ used in the deterministic smoothing algorithm {\em (right)}.\label{fig:path}}
\end{center}
\end{figure}

In Figure~\ref{fig:path} we plot the sequence of line search parameters $\gamma$ for the stochastic algorithm together with the values of the Lipschitz constant $L$ used in the deterministic smoothing algorithm, when solving problem~\eqref{eq:dspca} with $n=500$. We observe that both algorithms initially make longer steps, then slow down as they get closer to the optimum (where the leading eigenvalues are clustered).

%%%%%%%%%%% Switch appendix on/off %%%%%%%%%
%\iffalse

\section{Appendix}

In this Appendix, we recall several useful results related to the algorithm presented here. Subsection \ref{ss:lead-eig} summarizes the complexity of computing {\em one} leading eigenvector of a symmetric matrix (versus computing the entire spectrum).  In Subsection \ref{subsec:TechnicalDetailsSmoothedFk}, we prove a number of technical results concerning the function $\smoothedF_k$ and its components. In particular, we prove Theorem~\ref{th:lip-gap} linking the local Lipschitz constant of the gradient and the spectral gap. Finally, we show in Subsection \ref{ss:secularEqn} how the secular equation can be generalized to perturbations of higher rank, and we discuss extensions of our smoothing argument using GUE matrices.

\subsection{Computing one leading eigenvector of a symmetric matrix} \label{ss:lead-eig}
The complexity results detailed above heavily rely on the fact that extracting {\em one} leading eigenvector of a symmetric matrix $X\in\symm_n$ can be done by computing a few matrix vector products. This simple fact is easy to prove using the power method when the eigenvalues of $X$ are well separated, and Krylov subspace methods making full use of the matrix vector products converge even faster. However, the problem becomes more delicate when the spectrum of $X$ is clustered. The section that follows briefly summarizes how modern numerical methods produce eigenvalue decompositions in practice.

We start by recalling how packages such as LAPACK \cite{Ande99} form a full eigenvalue (or Schur) decomposition of a symmetric matrix $X\in\symm_n$. The algorithm is strikingly stable and, despite its $O(n^3)$ complexity, often competitive with more advanced techniques when the matrix $X$ is small. We then discuss the problem of approximating one leading eigenpair of $X$ using Krylov subspace methods with complexity growing as $O(n^2\log n)$ with the dimension (or less when the matrix is structured). In practice, we will see that the constants in these bounds differ significantly, with the cost of a full eigenvalue decompositions (and matrix exponentials) growing as $4n^3/3$ while computing one leading eigenpair has cost $cn^2$, with $c$ in the hundreds.

\subsubsection{Full eigenvalue decomposition.}
Full eigenvalue decompositions are computed by first reducing the matrix $X$ to symmetric tridiagonal form using Householder transformations, then diagonalizing the tridiagonal factor using iterative techniques such as the QR or divide and conquer methods for example (see \cite[Chap. 3]{Stew01} for an overview). The classical QR algorithm (see \cite[\S8.3]{Golu90}) implicitly relied on power iterations to compute the eigenvalues and eigenvectors of a symmetric tridiagonal matrix with complexity $O(n^3)$, however more recent methods such as the MRRR algorithm by \cite{Dhil03} solve this problem with complexity $O(n^2)$. Starting with the third version of LAPACK, the MRRR method is part of the default routine for diagonalizing a symmetric matrix and is implemented in the \texttt{STEGR} driver (see \cite{Dhil06}).

Overall, the complexity of forming a {\em full} Schur decomposition of a symmetric matrix $X\in\symm_n$ is then $4n^3/3$ flops for the Householder tridiagonalization, followed by $O(n^2)$ flops for the Schur decomposition of the tridiagonal matrix using the MRRR algorithm.

\subsubsection{Computing one leading eigenpair.} We now give a brief overview of the complexity of computing leading eigenpairs using Krylov subspace methods and we refer the reader to \cite[\S4.3]{Stew01}, \cite[\S8.3, \S9.1.1]{Golu90} or \cite{Saad92} for a more complete discussion. Successful termination of a {\em deterministic} power or Krylov method can never be guaranteed since in the extreme case where the starting vector is orthogonal to the leading eigenspace, the Krylov subspace contains no information about leading eigenpairs, so the results that follow are stochastic. \cite[Th.4.2]{Kucz92} show that, for any matrix $X\in\symm_n$ (including matrices with clustered spectrum), starting the algorithm at a random $u_1$ picked uniformly over the sphere means the Lanczos decomposition will produce a leading eigenpair with {\em relative} precision $\epsilon$, i.e. such that $|\lambda-\lambdamax|\leq \epsilon \lambdamax$, in
\[
k^\mathrm{Lan}\leq \frac{\log(n/\delta^2)}{4\sqrt{\epsilon}}
\]
iterations, with probability at least $1-\delta$. This is of course a highly conservative bound and in particular, the worst case matrices used to prove it vary with $k^\mathrm{Lan}$.

This means that computing one leading eigenpair of the matrix $X$ requires computing at most $k^\mathrm{Lan}$ matrix vector products (we can always restart the code in case of failure) plus $4nk^\mathrm{Lan}$ flops. When the matrix is dense, each matrix vector product costs $n^2$ flops, hence the total cost of computing one leading eigenpair of~$X$ is
\[
O\left(\frac{n^2\log(n/\delta^2)}{4\sqrt{\epsilon}}\right)
\]
flops. When the matrix is sparse, the cost of each matrix vector product is $O(s)$ instead of $O(n^2)$, where $s$ is the number of nonzero coefficients in $X$. Idem when the matrix $X$ has rank $r<n$ and an explicit factorization is known, in which case each matrix vector product costs $O(nr)$ which is the cost of two $n\times r$ matrix vector products, and the complexity of the Lanczos procedure decreases accordingly.

The numerical package ARPACK by \cite{Leho98} implements the Lanczos procedure with a reverse communication interface allowing the user to compute efficiently the matrix vector product $Xu_j$. However, it uses the implicitly shifted QR method instead of the more efficient MRRR algorithm to compute the Ritz pairs of the matrix $T_k\in\symm_k$.

\subsection{Technical results concerning $\smoothedF_k$ and its components}\label{subsec:TechnicalDetailsSmoothedFk}
\subsubsection{General remarks on rotational invariance}\label{ss:appRotaInvar}
We use repeatedly in the paper the fact that the type of smoothing we devised has some rotational invariance properties, which allows us to perform our computations on diagonal matrices without losing generality. We summarize the results we need and use in the following statement. 
\begin{lemma}\label{lemma:impactOfRotaInvariance}
Let $X$ be a deterministic matrix in $\symm_n$. $X$ is diagonalizable in an orthonormal basis and we write $X=\orthoMat_X\trsp D_X \orthoMat_X$, where $D_X$ is a diagonal matrix containing the eigenvalues of $X$ and $\orthoMat_X$ is a matrix of eigenvectors of $X$. Let $\{\canonicalRandomVector_i\}_{i=1}^k$ be $k$ i.i.d $\mathcal{N}(0,\idm_n)$ random vectors. Let $\nu$ be in $\reals$ and call 
$$
F_k(X)=\max_{1\leq i \leq k} \lambdamax(X+\nu \canonicalRandomVector_i\canonicalRandomVector_i\trsp)\;.
$$
Then 
$$
F_k(X)\stackrel{\mathcal{L}}{=} F_k(D_X)\;.
$$
Furthermore, if $\phi[F_k(X)]$ is an eigenvector associated with $F_k(X)$, we have 
$$
\phi[F_k(X)]\stackrel{\mathcal{L}}{=} \orthoMat_X\trsp \phi[F_k(D_X)]\;.
$$
\end{lemma}
\begin{proof}
	We observe that 
	$$
	X+\nu \canonicalRandomVector_i\canonicalRandomVector_i\trsp=\orthoMat_X\trsp \left[D_X+\nu (\orthoMat_X\canonicalRandomVector_i)(\orthoMat_X\canonicalRandomVector_i)\trsp\right] \orthoMat_X\;.
	$$
	Now it is a standard property of the normal distribution that if $\{\canonicalRandomVector_i\}_{i=1}^k$ are i.i.d $\mathcal{N}(0,\idm_n)$, then $\{\orthoMat_X\canonicalRandomVector_i\}_{i=1}^k$ are i.i.d $\mathcal{N}(0,\idm_n)$, for any (deterministic) orthonormal matrix $\orthoMat_X$. 
	The results we announced follow immediately. 
\end{proof}

\subsubsection{Existence of a density for $T$}\label{subsubsec:densityForT}
In Lemma \ref{lem:density}, we were interested in $T=\lambdamax(X+\mysigmasquared/n\canonicalRandomVector\canonicalRandomVector\trsp)-\lambdamax(X)$.
We prove Lemma \ref{lem:density} here, showing that $T$ has a density on $[0,\infty)$ when $\canonicalRandomVector \sim \mathcal{N}(0,\idm_n)$.

\begin{proof} [of Lemma \ref{lem:density}]
As usual, we call $\{\lambda_i\}_{i=1}^n$ the decreasingly ordered eigenvalues of $X$ and assume here that $\lambdamax(X)$ has multiplicity $l<n$ (if $l=n$ there is nothing to show, since then $X$ is proportional to $\idm_n$).
By rotational invariance of the standard Gaussian distribution, we can and do assume that $X$ is diagonal in what follows (see Lemma \ref{lemma:impactOfRotaInvariance} for details, if needed). As we have seen before, $T$ is therefore the only positive root of the equation
$$
0=s(T)=\frac{n}{\mysigmasquared}-\frac{\sum_{i=1}^l \coordCanonRandomVector_i^2}{T}-\sum_{i=l+1}^n\frac{\coordCanonRandomVector_i^2}{(\lambda_1-\lambda_i)+T},
$$
and note that $s(t)$ is increasing in $t$ when $t>0$. Hence, for any given $t>0$,
\begin{align*}
P(T\geq t)&=P(s(T)\geq s(t))=P(0\geq s(t))\\
&=P\left(\frac{\sum_{i=1}^l \coordCanonRandomVector_i^2}{t}+\sum_{i=l+1}^n\frac{\coordCanonRandomVector_i^2}{(\lambda_1-\lambda_i)+t}
\geq \frac{n}{\mysigmasquared}\right),\\
&=\int_{\frac{1}{\mysigmasquared}}^{\infty} p_t(u)du\triangleq I(t)\;,
\end{align*}
where $p_t$ is the density of the random variable
$$
Y_t=\frac{1}{n}\left(\frac{\sum_{i=1}^l \coordCanonRandomVector_i^2}{t}+\sum_{i=l+1}^n\frac{\coordCanonRandomVector_i^2}{(\lambda_1-\lambda_i)+t}\right).
$$
If the integral $I(t)$ can be differentiated under the integral sign, then we can differentiate $P(T\geq t)$ and we will have established the existence of a density for $T$ and hence for $\lambda_1+T$. Now, $p_t(x)$ is a very smooth function of both $t$ and $x$. Indeed, it is a convolution of $n-l$ densities that are smooth in $t$ and $x$. As a matter of fact, recall that if $X$ has density $p$ and $t>0$, $X/t$ has density $t p(t\cdot)$. Recall also that a random variable with $\chi^2_l$ distribution has density (see e.g \cite{MardiaKentBibby79}, p. 487) 
$$
p_l(x)=\frac{2^{-l/2}}{\Gamma(l/2)}x^{l/2-1}\exp(-x/2)\;  1_{x\in (0,\infty)}\; .
$$
So it is clear that for any $l$, any $x>0$, any $t>0$, and any $\alpha\geq 0$, $t\rightarrow (t+\alpha)p_l((t+\alpha)x)$ is $C^{\infty}$ in $t$. Applying this result in connection to \citep[Th.\,A.5.1]{Durr10}, we see that $Y_t$ has a density which is a smooth function of $t>0$. Indeed, it is $C^{\infty}$ on $(0,\infty)$. Moreover, it is easy to see that the conditions of \citep[Th.\,A.5.1]{Durr10} are satisfied for $p_t$, which guarantees that we can differentiate under the integral sign. This shows that for any $t>0$, the function $\pi$ such that $\pi(t)=P(T\geq t)$ is differentiable in $t$. It is also clear that $P(T=0)$ is 0, so this distribution has no atoms at 0. We conclude that $T$ has a density on~$[0,\infty)$.
\end{proof}

\subsubsection{Controlling the Hessian of $\lambdamax(X)$}\label{ss:Hessian}
Consider the map $F_0: \symm_n\rightarrow \mathbb{R}$ such that $F_0(X)=\lambdamax(X)$. We want to show that its gradient is Lipschitz continuous, when the largest eigenvalue of $X$ has multiplicity one and control the local Lipschitz constant. To do so, we compute 
$$
\gamma(X,Y)=\lim_{t\tendsto 0} \partial^2 F_0(X+tY)/\partial t^2\;,
$$ 
where $\norm{Y}_F=1$, and $Y$ is symmetric. It is standard that the local Lipschitz constant - with respect to Frobenius norm - of $\nabla F_0$ is 
$$
\localLipConstant{\nabla F_0(X)}=\sup_{Y \in \symm_n: \norm{Y}_F=1} \gamma(X,Y)\;.
$$

Let us call $\lambda_1>\lambda_2\geq \lambda_3\geq \ldots\geq\lambda_n$ the ordered eigenvalues of $X$. Very importantly we assume that $\lambda_1$ has multiplicity one. If not, it is easy to see that the function $\lambdamax(X)$ is continuous but not differentiable. We refer the reader to \citep{Kato95,Over95,Lewi02} for a more complete discussion.  Recall that in this situation Theorem \ref{th:lip-gap} stated that 
\begin{equation}\label{eq:valueLipConstantApp}
\localLipConstant{\nabla F_0(X)}=\frac{1}{\lambdamax(X)-\lambda_2(X)}.
\end{equation}
We now prove this statement.
% \begin{theorem}\label{thm:ControlOfHessian}
% Suppose $X$ is an $n\times n$ symmetric matrix with decreasingly ordered eigenvalues $\{\lambda_i\}_{i=1}^n$. Call $F_0(X)=\lambda_{\mathrm{max}}(X)$ and suppose that $\lambdamax(X)$ has multiplicity one. Let $Y$ be a symmetric matrix with $\norm{Y}_F=1$. Let us call
% $$
% g(X,Y)=\lim_{t\tendsto 0} \frac{\partial^2 F_0(X+tY)}{\partial t^2}\;.
% $$
% Then we have
% \begin{equation} \label{eq:valueLipConstantApp}
% \localLipConstant{\nabla F_0(X)}=\sup_{Y\in \symm_n,\norm{Y}_F=1} g(X,Y)=\frac{1}{2}\frac{1}{\lambda_1(X)-\lambda_2(X)}\;.
% \end{equation}
% \end{theorem}

\begin{proof} [of Theorem \ref{th:lip-gap}]
The strategy is to first exhibit a matrix $Y_c$ in $\symMatrices$ that will give us the right-hand side of Equation~\eqref{eq:valueLipConstantApp} as a lower bound. And then we will show that indeed this bound is the best one can do. We will rely heavily on the following classical result from the analytic perturbation theory of matrices. We can use \cite[p.81]{Kato95} or \citep{Lewi02} to get, for small $t$
$$
F_0(X+tY)=\lambdamax(X)+t \phi_1\trsp Y \phi_1+t^2 \sum_{j=2}^n \frac{1}{\lambda_1(X)-\lambda_j(X)} (\phi_1\trsp Y \phi_j)^2+\lo(t^2)\;,
$$
where $\phi_1$ is an eigenvector (of the matrix $X$) corresponding to the eigenvalue $\lambda_1$ and $\phi_j$ is an eigenvector (of $X$) corresponding to the eigenvalue $\lambda_j$. Here we have crucially used the fact that $\lambda_1(X)$ has multiplicity one. We conclude that 
\begin{equation}\label{eq:HessianInOneDirection}
\gamma(X,Y)=\lim_{t\tendsto 0} \frac{\partial^2 F_0(X+tY)}{\partial t^2}=2\sum_{j=2}^n \frac{1}{\lambda_1(X)-\lambda_j(X)} (\phi_1\trsp Y \phi_j)^2\;,
\end{equation}

\paragraph{\em Finding a lower bound for {$\localLipConstant{\nabla F_0(X)}$}}
Let $\orthoMat$ be an orthonormal matrix that transforms the canonical basis $(e_1,\ldots,e_n)$ into the orthonormal basis $(\phi_1,\ldots,\phi_n)$. In other words, $\orthoMat e_i=\phi_i$ and hence $\orthoMat\trsp \phi_i=e_i$. Let us call $P_0$ the matrix that exchanges $e_1$ and $e_2$ and send the other $e_j$'s to 0. In other words, the $2\times 2$ upper left block of $P_0$ is the matrix $\begin{pmatrix} 0&1\\1&0\end{pmatrix}$ and $P_0$ is zero everywhere else. Now call
$$
Y_c=\frac{1}{\sqrt{2}} \orthoMat P_0 \orthoMat\trsp \;.
$$
Note that $Y_c \in \symMatrices$. Since $\orthoMat\trsp \phi_i=e_i$, we see that $Y_c \phi_1 = \phi_2/\sqrt{2}$, $Y_c \phi_2 = \phi_1/\sqrt{2}$, and $Y_c \phi_j=0$ if $j>2$. Further, $\norm{Y_c}_F^2=\trace{Y_c\trsp Y_c}=\trace{Y_c^2}=\trace{\orthoMat P_0^2 \orthoMat\trsp}/2=\norm{P_0}_F^2/2=1$. Now, $\phi_1\trsp Y_c \phi_j=\delta_{2,j} \norm{\phi_1}^2/\sqrt{2}$. Hence,
$$
\gamma(X,Y_c)=\lim_{t\tendsto 0} \frac{\partial^2 F_0(X+tY_c)}{\partial t^2}=\frac{2}{2}\frac{1}{\lambda_1(X)-\lambda_2(X)}\;,
$$
and therefore,
$$
\localLipConstant{\nabla F_0(X)}\geq \frac{1}{\lambda_1(X)-\lambda_2(X)}\;.
$$

\paragraph{\em Finding an upper bound for {$\localLipConstant{\nabla F_0(X)}$}}
On the other hand, we clearly have, for $j\geq 2$, $0\leq 1/(\lambda_1(X)-\lambda_j(X))\leq 1/(\lambda_1(X)-\lambda_2(X))$. Therefore,
$$
\sum_{j=2}^n \frac{1}{\lambda_1(X)-\lambda_j(X)} (\phi_1\trsp Y \phi_j)^2\leq \frac{1}{\lambda_1(X)-\lambda_2(X)} \sum_{j=2}^n (\phi_1\trsp Y \phi_j)^2\;.
$$
Since $\{\phi_j\}_{j=1}^n$ form an orthonormal basis, and $Y$ is symmetric,
$$
\sum_{j=1}^n (\phi_1\trsp Y \phi_j)^2=\norm{Y\phi_1}_2^2\;.
$$
As a matter of fact $\phi_1\trsp Y \phi_j$ is just the coefficient of the vector $Y\trsp \phi_1=Y\phi_1$ in its representation in the basis of the $\phi_i$'s. We therefore have
$$
\sum_{j=2}^n \frac{1}{\lambda_1(X)-\lambda_j(X)} (\phi_1\trsp Y \phi_j)^2\leq\frac{1}{\lambda_1(X)-\lambda_2(X)}\left(\norm{Y\phi_1}_2^2-(\phi_1\trsp Y\phi_1)^2\right)\;.
$$
Let us call $\tilde{y}_{i,j}$ the $(i,j)$-th entry of the matrix that represents $Y$ in the basis of the $\phi_i$'s. Since $\norm{Y}_F^2=1$,
$$
\sum_{i,j}\tilde{y}_{i,j}^2=1\;.
$$
Using the symmetry of $Y$, we therefore see that
$$
2 \sum_{j=2}^n \tilde{y}_{1,j}^2 +\tilde{y}_{1,1}^2\leq 1\;.
$$
Now, $\norm{Y\phi_1}_2^2=\sum_{j=1}^n \tilde{y}_{1,j}^2$ and $(\phi_1\trsp Y\phi_1)^2=\tilde{y}_{1,1}^2$. Hence,
$$
\left(\norm{Y\phi_1}_2^2-(\phi_1\trsp Y\phi_1)^2\right)=\sum_{j=2}^n \tilde{y}_{1,j}^2\leq \frac{1-\tilde{y}_{1,1}^2}{2}\leq \frac{1}{2}\;.
$$
We conclude that
$$
\forall \, Y \in \symMatrices,\; \norm{Y}_F=1\;, \; \;\gamma(X,Y)\leq \frac{2}{2}\frac{1}{\lambda_1(X)-\lambda_2(X)}\;,
$$
and therefore
$$
\localLipConstant{\nabla F_0(X)}=\sup_{Y \in \symMatrices, \norm{Y}_F=1} \gamma(X,Y) \leq \frac{1}{\lambda_1(X)-\lambda_2(X)}\;.
$$
Since we have matching upper and lower bounds for $\localLipConstant{\nabla F_0(X)}$, we have established Theorem \ref{th:lip-gap}.
\end{proof}
%\redtext{WHAT FOLLOWS IS NEW; NEK 1/25/13; Alternative: just cite Chapter 7 of Schirotzek, ``Non-smooth analysis''}
\subsubsection{Differentials of maximum of several differentiable functions}\label{subsec:appDiffMaxSeveralFunctions}
We need the following elementary and well-known results at several points in the paper. We put them in this Appendix for the convenience of the reader. 
\begin{lemma}\label{lemma:diffSupSeveralFunctions}
Consider the function $\Psi_k=\max_{1\leq i \leq k} \psi_k$, where $\psi_1,\ldots,\psi_k$ are G\^ateaux-differentiable functions from $\mathcal{D}\subset \reals^d$ to $\reals$ and $k$ is an integer. Let $\textrm{int}(\mathcal{D})$ be the interior of $\mathcal{D}$. Let $x_0 \in \textrm{int}(\mathcal{D})$ be such that there exists $i_0 \in \{1,\ldots,k\}$ such that $\psi_{i_0}(x_0)>\psi_j(x_0)$ for all $j\neq i_0$. Then, $\Psi_k$ is G\^ateaux-differentiable at $x_0$ with 
$$
\nabla_G \Psi_k(x_0)= \nabla_G \psi_{i_0}(x_0)\;.
$$
Furthermore, when $\psi_j$'s are Fr\'echet-differentiable, so is $\Psi_k$ at $x_0$.
\end{lemma}
The proof shows that the result extends to higher order derivatives when they exist. \\
\begin{proof}
This is simply a restatement of the results of Proposition 7.2.7 in \citep{SchirotzekNonSmoothAnalysisBook07}, or \citep{HULemarechaL01} Theorem 4.4.2 and Corollary 4.4.4. 
We give the key idea and a proof of this easy fact for the sake of completeness. 

Indeed, let $I_j$ be the set of points $y$ such that $\psi_j(y)\geq \psi_{l}(y)$ for all $l\neq j$. We call $1_{I_j}$ the function taking value 1 on $I_j$ and 0 elsewhere. Let $N(x)$ be equal to $\card{j, 1\leq j \leq k: \psi_j(x)=\Psi_k(x)}$. Note that $N(x)=\sum_{j=1}^k 1_{I_j}(x)$. It is clear that $1\leq N(x)\leq k$. We also have 
$$
\Psi_k(x)=\frac{\sum_{j=1}^k \psi_j(x) 1_{I_j}(x)}{N(x)}\;.
$$
Under our assumptions on $x_0$, it is clear that $N(x_0)=1$. Furthermore, in a neighborhood $V(x_0)$ of $x_0$, we have $N(x)=1$ by continuity of the functions $\psi_j$'s. Of course, $V(x_0)$ is open, by definition of a neighborhood. It follows that for all $x$ in $V(x_0)$, we have $\Psi_k(x)=\psi_{i_0}(x)$. It now follows from the definition of G\^ateaux-differentiability that $\Psi_k$ is  G\^ateaux-differentiable at $x_0$ with the same G\^ateaux-differential as $\psi_{i_0}$. The result in the case of Fr\'echet differentiable functions $\psi_j$'s holds for the same reasons and is established by the same analysis.
\end{proof}

We have the following corollary. 
\begin{lemma}\label{lemma:localLipConstantLargestEigSeveralMatrices}
Suppose $X \in \symm_n$ and $v_1,\ldots,v_k$ are vectors in $\reals^n$. Denote by ${F}_k(X)=\max_{1\leq i \leq k} \lambdamax(X+v_{i}v_{i}\trsp)$. Suppose that $X$ and $\{v_i\}_{i=1}^k$ are such that there exists a unique $i_0$ such that 
$$
\lambdamax(X+v_{i_0}v_{i_0}\trsp)={F}_k(X)\;.
$$ 
Suppose further that the largest eigenvalue of $X+v_{i_0}v_{i_0}\trsp$ has multiplicity one. Then
$$
\localLipConstant{\nabla {F}_k(X)}=\frac{1}{\lambdamax(X+v_{i_0}v_{i_0}\trsp)-\lambda_2(X+v_{i_0}v_{i_0}\trsp)}\;.
$$
It follows that for $i_0=\argmax_{1\leq i \leq k}\lambdamax(X+v_{i}v_{i}\trsp)$, if ${F}_k(X)\neq \lambdamax(X)$,
$$
\localLipConstant{\nabla {F}_k(X)}\leq \frac{1}{\lambdamax(X+v_{i_0}v_{i_0}\trsp)-\lambdamax(X)}=\frac{1}{{F}_k(X)-\lambdamax(X)}\;.
$$
\end{lemma}
\begin{proof}
The proof of Lemma \ref{lemma:diffSupSeveralFunctions} shows that under our assumptions, $F_k$ coincides locally with $\lambdamax(X+v_{i_0}v_{i_0}\trsp)$. Hence the local Lipschitz constant of $\nabla {F}_k$ is the same as that of $X\mapsto \lambdamax(X+v_{i_0}v_{i_0}\trsp)$. One of our assumptions is that the largest eigenvalue of $X+v_{i_0}v_{i_0}\trsp$ has multiplicity 1. In that situation, Theorem \ref{th:lip-gap} guarantees that the local Lipschitz constant of $X\mapsto \lambdamax(X+v_{i_0}v_{i_0}\trsp)$ is 
$$
\frac{1}{\lambdamax(X+v_{i_0}v_{i_0}\trsp)-\lambda_2(X+v_{i_0}v_{i_0}\trsp)}\;.
$$
So we have established that 
$$
\localLipConstant{{F}_k(X)}=\frac{1}{\lambdamax(X+v_{i_0}v_{i_0}\trsp)-\lambda_2(X+v_{i_0}v_{i_0}\trsp)}\;.
$$
We now recall that by Cauchy's interlacing theorem (Theorem 4.3.4 in \citet{HornJohnson90}), $\lambda_2(X+v_{i_0}v_{i_0}\trsp)\leq \lambdamax(X)$. We can therefore conclude that 
$$
\localLipConstant{{F}_k(X)}=\frac{1}{\lambdamax(X+v_{i_0}v_{i_0}\trsp)-\lambdamax(X)}\;,
$$
since we have assumed that $\lambdamax(X)\neq \lambdamax(X+v_{i_0}v_{i_0}\trsp)={F}_k(X)$.
\end{proof}

\subsubsection{Interchanging expectation and differentiation for $\smoothedF_k$}\label{ss:expAndDiffInterchange}
\begin{lemma}\label{lemma:interchangeDiffAndEexpect}
We can interchange expectation and differentiation for $\smoothedF_k$ so that
$$
\nabla \smoothedF_k(X)=\Exp{\phi_{i_0}\phi_{i_0}\trsp}
$$
using the notation of Lemma \ref{eq:grad-var}.
\end{lemma}
\newcommand{\spectralNorm}[1]{\norm{#1}_2} % AA: \|.\|_2 is standard in optim
\begin{proof}
$\smoothedF_k$ is convex as an average of convex functions. To show that it is differentiable, it is therefore enough to show that it is G\^ateaux-differentiable (see \cite{HULemarechaL01}, Corollary D.2.1.4). 	
Let $X_0$ be given. We use the notation 
$$
\smoothedF_k(X_0)=\int F_k(X_0;\canonicalRandomVector)\mu(d\canonicalRandomVector)
$$ 
to make things simpler in this proof. Of course, $F_k(X_0;z)=\max_{1\leq i \leq k} \lambdamax(X_0+\frac{\eps}{n}\canonicalRandomVector_i\canonicalRandomVector_i\trsp)$. (Compared to the main text, we have now made the dependence on $\canonicalRandomVector$ explicit as it is needed below to address a potential measure theoretic problem.) $\mu(d\canonicalRandomVector)$ is just the joint distribution of $\canonicalRandomVector_i$'s, for $i=1,\ldots,k$. To make notations simple in this proof, we use $z$ to denote $(z_i)_{i=1}^k$.\\
We know that $F_k(X_0;\canonicalRandomVector)$ has a subdifferential for all $X_0$ and all $\canonicalRandomVector$, since it is the maximum of $k$ functions with a subdifferential (see \cite{HULemarechaL01}, Theorem D.4.4.2). The spectral norm of the elements of this subdifferential is bounded by 1, since they are a convex combination of matrices of spectral norm at most 1 (see \cite{HULemarechaL01}, Theorem D.4.4.2 and Equation (5.1.3) p. 195 in that book, which characterizes the subdifferential of the largest eigenvalue mapping of a symmetric matrix). \\
Suppose $Y_0$ is a fixed matrix, with $\spectralNorm{Y_0}=1$ without loss of generality, where $\spectralNorm{Y}$ is the spectral norm of the symmetric matrix $Y$. By the mean value theorem for functions with a subdifferential (see \cite{HULemarechaL01}, Theorem D.2.3.4) we have, if we call $X_{0,t}=X_0+tY_0$
$$
F_k(X_0+tY_0;\canonicalRandomVector)-F_k(X_0;\canonicalRandomVector)=t\int_0^1 <\partial F_k(X_{0,tu};\canonicalRandomVector),Y_0> du\;,
$$
where $\partial F_k(X_{0,tu};\canonicalRandomVector)$ is any choice of subgradient of $F_k(X_{0,tu};\canonicalRandomVector)$ and $<A,B>=\trace{A\trsp B}$ for the symmetric matrices we are working with. (We use the same notation in this proof for subgradients and subdifferentials since it does not create confusion.)\\
Because $\spectralNorm{\partial F_k(X_{0,tu})}\leq 1$, we can apply Fubini's theorem to get 
$$
\frac{\smoothedF_k(X_0+tY_0)-\smoothedF_k(X_0)}{t}=\int_0^1 du \int <\partial F_k(X_{0,tu};\canonicalRandomVector),Y_0> \mu(d\canonicalRandomVector)\;.
$$
% We have shown (see Corollary \ref{cor:diff})  that for any fixed $tu \in \mathbb{R}$, there exists a set $E_{tu}$ such that $\partial F_k(X_{0,tu};\canonicalRandomVector)$ is single valued (i.e $F_k(X_{0,tu};\canonicalRandomVector)$ is differentiable) with $\mu(E_{tu})=1$. However, we cannot conclude a priori that $\mu(\bigcap_{u}E_{tu})=1$, since the index set is not denumerable. This creates a potential problem. \\
Now let $\eta>0$ be given and let 
\begin{align*}
\mathcal{A}_\eta&=\{\canonicalRandomVector: \; \forall i\;, 1\leq i \leq k\;, \lambdamax(X+\frac{\eps}{n}\canonicalRandomVector_i\canonicalRandomVector_i\trsp)-\lambdamax(X)>\eta\}\;,\\
\mathcal{B}_\eta&=\{\canonicalRandomVector: \; \exists i_0: \forall j\neq i_0\;, 1\leq j \leq k\;, \lambdamax(X+\frac{\eps}{n}\canonicalRandomVector_{i_0}\canonicalRandomVector_{i_0}\trsp)\geq \lambdamax(X+\frac{\eps}{n}\canonicalRandomVector_j\canonicalRandomVector_j\trsp)+\eta\}\;,\\
\mathcal{E}_\eta&=\mathcal{A}_\eta \cap \mathcal{B}_\eta\;.
\end{align*}
By continuity of the maps involved, $\mathcal{E}_\eta$ is clearly measurable with respect to Lebesgue measure and therefore~$\mu$. Equation \eqref{eq:gap-lb} in Proposition \ref{prop:impactLowRankUpdateOnTopEigenvalue} implies that $\lim_{\eta\tendsto 0} \mu(\mathcal{A}_\eta)=1$. Lemma \ref{lem:density} also implies that $\lim_{\eta\tendsto 0} \mu(\mathcal{B}_\eta)=1$. We conclude that $\lim_{\eta\tendsto 0} \mu({\mathcal E}_\eta)=1$.

This implies that  
$$
\left|\int <\partial F_k(X_{0,tu};\canonicalRandomVector),Y_0> \mu(d\canonicalRandomVector)-\int_{\mathcal{E}_\eta} <\partial F_k(X_{0,tu};\canonicalRandomVector),Y_0> \mu(d\canonicalRandomVector)\right|\leq n \mu(\mathcal{E}^c_\eta)\tendsto 0 \text{ as } \eta\tendsto 0\;,
$$
since the absolute value of the integrand is bounded by $n$. For the same reasons, as $\eta\tendsto 0$,
$$
\left|\int_{0}^1du \int <\partial F_k(X_{0,tu};\canonicalRandomVector),Y_0> \mu(d\canonicalRandomVector)-\int_0^1 du\int_{\mathcal{E}_\eta} <\partial F_k(X_{0,tu};\canonicalRandomVector),Y_0> \mu(d\canonicalRandomVector)\right|\leq n \mu(\mathcal{E}^c_\eta)\tendsto 0.
$$
When $|t|<\eta/4$, it is clear that for all $\canonicalRandomVector \in \mathcal{A}_\eta$, and all $u\in [0,1]$, we have 
$$
\forall i\;, 1\leq i \leq k\;, \lambdamax(X+tuY_0 +\frac{\eps}{n}\canonicalRandomVector_i\canonicalRandomVector_i\trsp)-\lambdamax(X+tuY_0)>\eta/2\;.
$$
As matter of fact, we have 
\begin{align*}
\lambdamax(X+tuY_0 +\frac{\eps}{n}\canonicalRandomVector_i\canonicalRandomVector_i\trsp)-\lambdamax(X+tuY_0)&=
\lambdamax(X+tuY_0 +\frac{\eps}{n}\canonicalRandomVector_i\canonicalRandomVector_i\trsp)-
\lambdamax(X+ \frac{\eps}{n}\canonicalRandomVector_i\canonicalRandomVector_i\trsp)\\
&+
\lambdamax(X+ \frac{\eps}{n}\canonicalRandomVector_i\canonicalRandomVector_i\trsp)-
\lambdamax(X)\\
&+\lambdamax(X)-
\lambdamax(X+tuY_0)
\end{align*}
By Weyl's inequality, $|\lambdamax(X+tuY_0 +\frac{\eps}{n}\canonicalRandomVector_i\canonicalRandomVector_i\trsp)-\lambdamax(X+ \frac{\eps}{n}\canonicalRandomVector_i\canonicalRandomVector_i\trsp)|\leq |t|u \norm{Y_0}\leq \eta/4$ and $|\lambdamax(X)-\lambdamax(X+tuY_0)|\leq \eta/4$ for the same reason. By the same reasoning, if $z\in \mathcal{B}_\eta$, when $|t|<\eta/4$, for all $u\in [0,1]$,
$$
\lambdamax(X+tuY_0 +\frac{\eps}{n}\canonicalRandomVector_{i_0}\canonicalRandomVector_{i_0}\trsp)\geq \lambdamax(X+tuY_0+ \frac{\eps}{n}\canonicalRandomVector_j\canonicalRandomVector_j\trsp)+\frac{\eta}{2}\;.
$$

This shows that for $\canonicalRandomVector \in \mathcal{E}_\eta$, when $|t|<\eta/4$, $F_k(X_{0,tu};\canonicalRandomVector)$ is differentiable for all $u$ and so its subdifferential is reduced to a singleton. Therefore, 
$$
\forall \canonicalRandomVector \in \mathcal{E}_\eta \,, 0\leq u \leq 1,\;  \text{ and } |t|<\eta/4 \;, \; \;\partial F_k(X_{0,tu};\canonicalRandomVector)=\nabla F_k(X_{0,tu};\canonicalRandomVector)\;.
$$
We know in fact that under the aforementioned conditions $F_k(X_{0,tu};\canonicalRandomVector)$ is twice differentiable as a function of $tu$ (see \cite{Kato95}, pp.80-81) and therefore the gradient $\nabla F_k(X_{0,tu};\canonicalRandomVector)$ is continuous (as a function of $tu$). So we have, pointwise in $\canonicalRandomVector \in \mathcal{E}_\eta$ and $u\in [0,1]$,
$$
\lim_{t\tendsto 0} \nabla F_k(X_{0,tu};\canonicalRandomVector)=\nabla F_k(X_{0};\canonicalRandomVector)\;.
$$
Using the fact that $\nabla F_k(X_{0,tu};\canonicalRandomVector)$ is bounded (it is a rank 1 matrix of norm 1), we can use the dominated convergence theorem and Fubini's theorem to conclude that 
%for all $u \in [0,1]$,
% $$
% \lim_{t\tendsto 0} \int_{\mathcal{E_\eta}} <\partial F_k(X_{0,tu};\canonicalRandomVector),Y_0> \mu(d\canonicalRandomVector)=
% \int_{\mathcal{E_\eta}} <\nabla F_k(X_0;\canonicalRandomVector),Y_0> \mu(d\canonicalRandomVector)\;,
% $$
% and 
\begin{align*}
\lim_{t\tendsto 0} \int_0^1 du \int_{\mathcal{E_\eta}} <\partial F_k(X_{0,tu};\canonicalRandomVector),Y_0> \mu(d\canonicalRandomVector)&=
\lim_{t\tendsto 0} \int_0^1 du \int_{\mathcal{E_\eta}} <\nabla F_k(X_{0,tu};\canonicalRandomVector),Y_0> \mu(d\canonicalRandomVector)\\
&=
\int_{\mathcal{E_\eta}} <\nabla F_k(X_0;\canonicalRandomVector),Y_0> \mu(d\canonicalRandomVector)\;.
\end{align*}
Finally, because $<\partial F_k(X_0;\canonicalRandomVector),Y_0>$ is bounded and because $\mathcal{E}_\eta$ is a decreasing family of sets, we see that 
$$
\lim_{\eta\tendsto 0}\int_{\mathcal{E_\eta}} <\nabla F_k(X_0;\canonicalRandomVector),Y_0> \mu(d\canonicalRandomVector)=
\int <\partial F_k(X_0;\canonicalRandomVector),Y_0> \mu(d\canonicalRandomVector)\;,
$$
by the dominated convergence theorem. Naturally, the previous equality is true for any choice of subgradients on the set of ($\mu$-)measure 0 where 
the subdifferential $F_k(X_0;\canonicalRandomVector)$ is not reduced to a singleton.
Furthermore,
$$
\lim_{t\tendsto 0}\frac{\smoothedF_k(X_0+tY_0)-\smoothedF_k(X_0)}{t}=<\int \partial F_k(X_0;\canonicalRandomVector)\mu(d\canonicalRandomVector),Y_0>\;.
$$
So we conclude that $\smoothedF_k$ is G\^ateaux-differentiable at $X_0$ and hence differentiable, since $\smoothedF_k$ is convex. The previous expression is valid for any subgradient of $F_k(X_0;\canonicalRandomVector)$. Since with probability 1 the subdifferential is a singleton, we can also write 
$$
\nabla \smoothedF_k(X_0)=\int \nabla(F_k(X_0;\canonicalRandomVector)) \mu(d\canonicalRandomVector)
$$
In the particular situation we are considering, this can also be re-written as 
$$
\nabla \smoothedF_k(X_0)=\Exp{\phi_{i_0}\phi_{i_0}\trsp}\;.
$$
which is the desired result.
\end{proof}%}

We now show that the gradient is diagonalizable in the same basis as $X$.

\begin{lemma}\label{lem:sim-diag}
The matrix $\nabla \smoothedF_k(X)$ is diagonalizable in the same basis as $X$. In particular, when $X$ is diagonal, so is $\nabla \smoothedF_k$.
\end{lemma}
\begin{proof}
Call $\lambda_i(X)$ the eigenvalues of $X$ in decreasing order. As above,  $\canonicalRandomVector_i\iid {\mathcal N}(0,\id_n)$ implies that no $\canonicalRandomVector_i$ is an eigenvector of $X$ with probability one.  We call $l_{1,i}=\lambdamax(X+\frac{\eps}{n}z_i z_i\trsp)$ and $\phi_i$ the corresponding eigenvector. With probability 1, $\phi_{i}$ is uniquely defined, up to sign, since $l_{1,i}$ has multiplicity 1 with probability 1.

\textbf{(i) $X$ diagonal.} We first focus on the case where $X$ is diagonal. Our strategy is to show that the off-diagonal entries of $\phi_{i_0}\phi_{i_0}\trsp$ have a (marginal) distribution that is symmetric around 0. 

%We know that $\nabla g(X,\canonicalRandomVector)=\phi(X,\canonicalRandomVector)\phi(X,\canonicalRandomVector)^T$ with
In what follows we use the notation $\canonicalRandomVector_i(j)$ to denote the $j$-th coordinate of the vector $\canonicalRandomVector_i$. 
It is well-known (\citep[see][\S8.5.3]{Golu90}, Theorem 8.5.3) that when $X$ is diagonal, the $j$-th coordinate of $\phi_i$ is given by
\begin{equation}\label{eq:coordsEigenVecRankOnePerturb}
\phi_i(j)=c\frac{\canonicalRandomVector_i(j)}{l_{1,i}-\lambda_j},
\end{equation}
where  $c>0$ is a normalizing factor. Recall that $l_{1,i}$ is the largest root of $\chi(\lambda)=0$, where
\[%\label{eq:SecularEquation}
\chi(\lambda)=1+\frac{n}{\mysigmasquared} \frac{\sum_{j=1}^l[\canonicalRandomVector_i(j)]^2}{\lambda_1(X)-\lambda}+\frac{n}{\mysigmasquared} \sum_{j=l+1}^n \frac{[\canonicalRandomVector_i(j)]^2}{\lambda_j(X)-\lambda}\;.
\]
This equation shows in particular that $l_{1,i}$'s depend on $\canonicalRandomVector_i$'s only through the absolute values of the coordinates of these vectors.
Let us now pick $j_0$, an integer such that $1\leq j_0 \leq n$. Suppose that we change the $j_0$-th coordinates of the vectors $\canonicalRandomVector_i$'s to their opposites. Call $\tilde{l}_{1,i}$ and $\tilde{\phi}_i$ the corresponding eigenvalue and eigenvectors. As we have just seen, 
$$
\tilde{l}_{1,i}=l_{1,i}\;, \forall i\;.
$$
In particular, $i_0$ is unaffected by this sign change operation.\\ 
On the other hand, 
$$
\forall i\;, \; \; 
\left[\tilde{\phi}_i\tilde{\phi}_i\trsp\right](l,m)=
\left\{
\begin{array}{cl}
-\left[\phi_i\phi_i\trsp\right](l,m)& \text{ if } l\neq m \text{ and } m=j_0 \text{ or } l=j_0\;,\\
\left[\phi_i\phi_i\trsp\right](l,m)& \text{ otherwise}\;.
\end{array}
\right.
$$
However, since $\canonicalRandomVector_i$'s have a symmetric distribution, their distribution is unaffected by a change of sign to one of the coefficients. So it is clear that 
$$
\forall i\;, \; \; \tilde{\phi}_i\tilde{\phi}_i\trsp\equalInLaw \phi_i\phi_i\trsp\;.
$$
So for all $i$, the distribution of the off-diagonal entries of the matrix $\phi_i\phi_i\trsp$ is symmetric around 0, since it is equal in law to its opposite. (We have just shown it for the off-diagonal entries for the $j_0$-th row and columns, but since there was nothing special about $j_0$, it is true for all the off-diagonal entries.) Furthermore, since the value of $i_0$ is unaffected by the sign change operation we discussed, we have shown that the off-diagonal entries of the matrix $\phi_{i_0}\phi_{i_0}\trsp$ have a symmetric distribution. Since $\phi_{i_0}\trsp \phi_{i_0}=1$, the entries of the matrix $\phi_{i_0}\phi_{i_0}\trsp$ are bounded and therefore have a mean. This mean must be zero for the off-diagonal entries since they have a symmetric distribution. So we have shown that $\Exp{\phi_{i_0}\phi_{i_0}\trsp}$ is a diagonal matrix when $X$ is diagonal.

\textbf{(ii) $X$ not diagonal.} When $X$ is not diagonal, we simply diagonalize $X$ into $X=\orthoMat_X\trsp D_X \orthoMat_X$ and use rotational invariance of the distribution of the $\canonicalRandomVector_i$'s to see that 
$$
\left[\phi_{i_0}(X)\phi_{i_0}(X)\trsp\right]\equalInLaw \orthoMat_X\trsp \left[\phi_{i_0}(D_X)\phi_{i_0}(D_X)\trsp\right] \orthoMat_X\;,
$$
where by a slight abuse of notation we have denoted by $\phi_{i_0}(X)$ an eigenvector associated with $\max_{1\leq j \leq k} l_{1,j}$. 
Since we have already seen that $\Exp{\phi_{i_0}(D_X)\phi_{i_0}(D_X)\trsp}$ is diagonal, we have shown that $\Exp{\phi_{i_0}(X)\phi_{i_0}(X)\trsp}$ is diagonal in the basis that diagonalizes $X$. 
\end{proof}

\subsection{On the secular equation and higher-order perturbations}\label{ss:secularEqn}
We give an elementary proof of the validity of the secular equation, which avoids matrix representations. Though simple and possibly well-known, the advantage of our derivation is that it extends easily to higher rank perturbation.
 More precisely, let us consider the matrix
 \begin{equation}
 M_1=M+U\;,
 \end{equation}
where $U$ is a symmetric matrix. We assume without loss of generality that $M$ is diagonal. We write $U=\sum_{j=1}^k \canonicalDetVector_j \canonicalDetVector_j^{T}$. We do not require the $\canonicalDetVector_j$ to be orthogonal and they could also be complex valued in what follows.

Let us call $\lambda_1\geq \lambda_2\geq \ldots \geq \lambda_n$ the eigenvalues of $M$. Our aim is to compute the characteristic polynomial of $M_1$ and relate it to that of $M$. We call
 \begin{align*}
 P_{M_1}(\lambda)&=\det(M_1-\lambda \id_n)\;,\\
 P_{M}(\lambda)&=\det (M-\lambda \id_n)\;,\\
 M_{\lambda}&=M-\lambda \id_n\;.
 \end{align*}
Assuming for a moment that $\lambda$ is not an eigenvalue of $M$, we clearly have $M_1-\lambda \id_n =M_{\lambda}(\id_n+M_{\lambda}^{-1}U)$. We call $G(\lambda)$ the $k \times k$ matrix with $(i,j)$ entry $\canonicalDetVector_j^{T} M_{\lambda}^{-1}\canonicalDetVector_{i}$.

We have
 $$
 P_{M_1}(\lambda)=\det(M_{\lambda})\det(\id_n+M_{\lambda}^{-1}U)=P_{M}(\lambda)\det(\id_k+G(\lambda))\;,
 $$
since $\det(\id_n+AB)=\det(\id_k+BA)$ for rectangular matrices $A$ and $B$ whenever $AB$ is $n\times n$ and $BA$ is $k\times k$.
The previous formula can be used to study the eigenvalues of finite rank perturbations of $M$, since they are the zeros of the characteristic polynomial $P_{M_1}$.

Let us focus on the case where $U$ has rank one, that is $U=\canonicalDetVector\canonicalDetVector\trsp$. Since we assume wlog that $M$ is diagonal, we have, when $k=1$,
$$
\det(\id_k+G(\lambda))=\det(1+\canonicalDetVector\trsp M_{\lambda}^{-1}\canonicalDetVector)=1+\sum_{i=1}^n \frac{\coordCanonDetVector_i^2}{\lambda_i-\lambda}\;.
$$
We therefore get, when $\lambda$ is not an eigenvalue of $M$,
\begin{equation}\label{eq:CharacPolynomial}
P_{M_1}(\lambda)=\left[\prod_{i=1}^n (\lambda_i-\lambda) \right]\left(1+\sum_{i=1}^n \frac{\coordCanonDetVector_i^2}{\lambda_i-\lambda}\right)\;,
\end{equation}
from which the secular equation follows. From Equation (\ref{eq:CharacPolynomial}), it is also clear that if $\lambda_i$ is an eigenvalue of $M$ with multiplicity $m>1$, $\lambda_i$ is also an eigenvalue of $M_1$ with multiplicity $m-1$.

\subsection{GUE smoothing}\label{s:speculative-math}

In this section, we discuss possible extensions of the stochastic regularization techniques, their efficiency and regularity. We have chosen to analyze the rank one perturbation scheme - and slight variants of it - because of its numerical efficiency and mathematical simplicity. However, many other random smoothing algorithms are possible and modern random matrix theory offers tools to understand their properties. We expect that some of them will lead to better worst case bounds than the order $n$ bound on the Lipschitz constant of the gradient for the rank one Gaussian perturbations we have considered here.

A case in point is the following. Consider a matrix $U$ from the Gaussian Unitary Ensemble (GUE). Matrices from $GUE$ are Hermitian random matrices with complex Gaussian entries, i.i.d $\complexNormal{0}{1}$ above the diagonal and i.i.d $\centeredNormal{1}$ on the diagonal. Recall that if $\coordCanonRandomVector_{\mathbb{C}}$ is $\complexNormal{0}{1}$, $\coordCanonRandomVector_{\mathbb{C}}=(\coordCanonRandomVector_1+i\coordCanonRandomVector_2)/\sqrt{2}$, where $\coordCanonRandomVector_1$ and $\coordCanonRandomVector_2$ are independent with distribution $\centeredNormal{1}$.

In what follows, $X$ is a deterministic matrix and $U$ is a random GUE matrix. We assume, without loss of generality, that the largest eigenvalue of $X$ is bounded, for instance $\lambda_{\max(X)}=1$ (if not, we can always shift $X$ by a multiple of $\idm_n$, which takes care of the problem).

A natural smoothing of $\lambda_{\max}(X)$ is $\smoothedF_{GUE}(X)=\Exp{\lambda_{\max}(X+(\mysigma/\sqrt{n})U)}$, where $U$ is a GUE matrix. This type of matrices belong to the so-called ``deformed GUE". \cite{johanssonGumbelToTW} is an important paper in this area and contains a result, Theorem 1.12, that is not exactly suited to our problem but quite close, perhaps despite the appearances. Before we proceed, we note that showing that $\smoothedF_{GUE}(X)$ is an $\mysigma$-approximation of $\lambdamax(X)$ is immediate from standard results on $GUE$ matrices (see \cite{Trotter84}, \cite{DavidsonSzarek01}).

In a nutshell, random matrix theory indicates that $\lambda_{\max}(X+({\mysigma}/{\sqrt{n}})U)$ undergoes a phase transition as  $\mysigma$ changes when $X$ is not a multiple of $\idm_n$. If $\mysigma$ is sufficiently large (more details follow), the behavior of $\lambda_{\max}(X+({\mysigma}/{\sqrt{n}})U)$ is driven by the GUE component and the spacing between the two largest eigenvalues is of order $n^{-2/3}$. On the other hand, if $\mysigma$ is not large enough, we remain essentially in a perturbative regime and the spacing between the two largest eigenvalues is larger than $n^{-2/3}$. A very detailed study of the phase transition should be possible, too. However, all these results are asymptotic. Non-asymptotic results could be obtained (the machinery to obtain results such as Johansson's is non-asymptotic) but would be hard to interpret and exploit. We therefore keep this discussion at an informal level.

Smoothing by a GUE matrix should give a worst case bound on $\localLipConstant{\nabla \smoothedF_{GUE}}$ of order $n^{2/3}$, which is better than the worst case bound of $n$ we have when we smooth with rank one matrices (but requires generating $O(n^2)$ random numbers instead of $O(n)$). GUE smoothing might therefore improve the performance of the algorithm since the cost of generating these random variables is typically dominated by the cost of computing a leading eigenvector of the perturbed matrix.

Let us give a bit more quantitative details. Based on Johansson's work and the solution to a similar problem in a different context (\cite{nekGenCov}), it is clear that the condition for the spacings to be of order $n^{-2/3}$ is the following (this result might be available in the literature but we have not found a reference). Call $\Lambda_n$ the spectral distribution of the $n\times n$ matrix $X_n$, i.e the probability distribution that puts mass $1/n$ at each of the $n$ eigenvalues of $X_n$. Call $w_c$ the solution in $(\lambda_{\max}(X_n),\infty)$ of
$$
\int \frac{d\Lambda_n(t)}{(w_c-t)^2}=\frac{1}{\mysigma^2}\;.
$$
Call $\goodClass$ the class of matrices for which
$$
\liminf_{n\rightarrow \infty} \left[w_c-\lambdamax(X_n)\right]>0\;.
$$

Then, looking carefully at Johansson's and El Karoui's work, it should be possible to show that:
if the sequence of matrices $X_n$ is in $\goodClass$, then, if $X_n(\mysigma)=X_n+\mysigma/\sqrt{n}U$,
$$
n^{2/3} \frac{\lambda_{\max}(X_n(\mysigma))-\alpha_n}{\beta_n} \weakCV \twc\;,
$$
where
$$
\alpha_n=w_c+\mysigma^2 \int \frac{d\Lambda_n(t)}{w_c-t}
\quad \mbox{and} \quad
\beta_n=\mysigma^2 \left(\int \frac{d\Lambda_n(t)}{(w_c-t)^3}\right)^{1/3}
$$
and $\twc$ is the Tracy-Widom distribution appearing in the study of GUE \citep[see][]{tw94}. The same is true for the joint distribution of the $k$ largest eigenvalues, where $k$ is a fixed integer, and $\twc$ is replaced by the corresponding limiting joint distribution for the $k$ largest eigenvalues of a GUE matrix.

When the matrix $X_n$ is not in $\goodClass$, then the top two eigenvalues should have spacing greater than $n^{-2/3}$. We expect that if $X_n$ has some sufficiently separated eigenvalues with multiplicity higher than one, the spacings there are at least $n^{-1/2}$, by analogy with \cite{CapitaineLargestEigenvaluesAOP09} and \cite{bbap}. To quantify what ``sufficiently separated" means, we could suppose that $X_n$ is a completion of a $(n-k_0) \times (n-k_0)$ matrix $X_{n-k_0,0}$ which is in $\goodClass$, to which we add $k_0$ eigenvalues $\lambda_{\max}(X_n)$, all equal and greater than $\lambda_{\max}(X_{n-k_0,0})$, with $\lambda_{\max}(X_n)$ greater than and bounded away from $w_c(X_{n-k_0,0})$. Calling $\Lambda_{n-k_0,0}$ the spectral distribution of $X_{n-k_0,0}$, we should have
$$
n^{1/2} \frac{\lambda_{\max}(X_n(\mysigma))-\alphaTilde_n}{\betaTilde_n} \weakCV \lambda_{\max}\left(\mathrm{GUE}_{k_0\times k_0}\right)\;,
$$
where $\alphaTilde_n=\lambda_{\max}(X_n)+\mysigma^2 \int \frac{d\Lambda_{n-k_0,0}(t)}{\lambda_{\max}(X_n)-t}$ and
$\betaTilde_n= \mysigma \left(1-\mysigma^2 \int \frac{d\Lambda_{n-k_0,0}(t)}{(\lambda_{\max}(X_n)-t)^2}\right)^{1/2}$.

The same is true for the $k_0$ largest eigenvalues of $X_n(\mysigma)$ and $\lambda_{\max}(\mathrm{GUE}_{k_0\times k_0})$ is replaced by the corresponding joint distribution for the $k_0\times k_0$ GUE.

In light of the integrability problems we had in the rank one perturbation case for the inverse spectral gap $1/(l_1(X_n(\mysigma))-l_2(X_n(\mysigma)))$, it is natural to ask whether such problems would arise with a GUE smoothing. For this informal discussion, we limit ourselves to answering this question for GUE (and not deformed GUE).
We recall that the joint density of the eigenvalues $\{l_{i,GUE}\}_{i=1}^n$ of a $n\times n$ GUE matrix is
$$
C \exp(-\sum_{i=1}^n l_{i,GUE}^2/2) \prod_{1\leq i<j\leq n} |l_{i,GUE}-l_{j,GUE}|^2 \;,
$$
where $C$ is a normalizing constant. So we see immediately that $1/(l_{1,GUE}-l_{2,GUE})$ is integrable in the GUE setting.
(The formula above is often stated for the unordered eigenvalues of a GUE matrix. The functional form of the density is unchanged by ordering, because of the symmetry. The domain of definition and the constant change when considering ordered eigenvalues, but this has no bearing on the question of integrability.)

The smoothing could also be done by a matrix from the Gaussian Orthogonal Ensemble (GOE), where the entries above the diagonal are i.i.d ${\mathcal N}(0,1)$ and the entries on the diagonal are i.i.d ${\mathcal N}(0,2)$. We do not know of a result corresponding to Johansson's in that case, though we would expect that the behavior of the top eigenvalues is the same as described above, with $\twc$ replace by $\tw_1$, the Tracy-Widom distribution appearing in the study of GOE.  From an algorithmic point of view, the two methods should therefore be equivalent.

%%%%%%%%%%%%% Switch appendix
%\fi

\section*{Acknowledgments}  Alexandre d'Aspremont would like to acknowledge partial support from NSF grants SES-0835550 (CDI), CMMI-0844795 (CAREER), CMMI-0968842, a starting grant from the European Research Council (project SIPA), a Peek junior faculty fellowship, a Howard B. Wentz Jr. award and a gift from Google. Noureddine El Karoui acknowledges support from an Alfred P. Sloan research Fellowship and NSF grant DMS-0847647 (CAREER).

\small{\bibliographystyle{plainnat}\bibsep 1ex
%\bibliography{biblioNEK}}
\bibliography{/Users/aspremon/Dropbox/Research/Biblio/MainPerso}}

\end{document}